\documentclass[12pt,a4paper,oneside]{article}
\usepackage[T1]{fontenc}
\usepackage[utf8]{inputenc}
\usepackage{amsmath,amsthm,amssymb,bm,graphicx,hyperref,mathrsfs}
\newtheorem{definition}{Definition}
\newtheorem{theorem}{Theorem}
\newtheorem{lemma}{Lemma}
\newtheorem{proposition}{Proposition}
\newtheorem{corollary}{Corollary}
\newtheorem{remark}{Remark}
\newcommand{\keywords}[1]{\par\noindent\textbf{Keywords:} #1}
\newcommand{\ElDiag}{\operatorname{ElDiag}}
\newcommand{\Ord}{\operatorname{Ord}}
\newcommand{\Trans}{\operatorname{Trans}}
\newcommand{\Mstar}{\mathcal{M}^{*}}
\newcommand{\PowTwo}{\operatorname{Pow}_2}
\newcommand{\Bit}{\operatorname{Bit}}

\newcommand{\B}{\mathsf{B}}

\title{Characterizing Fragments of Collection in Set Theory by Model-Theoretic Properties}
\author{Junhong Chen}
\date{\today}

\begin{document}

\bibliographystyle{plain}

\maketitle

\begin{abstract}
We compare the model theory of the weak set theory $\mathsf{DB}_0$ with
that of the algebraic theory $\mathsf{PA}^-$.  Every model of
$\mathsf{DB}_0$ has a proper taller* end extension with an exact
transitive cover, paralleling the corresponding end-extension result for
$\mathsf{PA}^-$.  More substantially, we prove that $\mathsf{DB}_0$ and
$\mathsf{PA}^-$ are mutually interpretable.  The new direction interprets
$\mathsf{DB}_0$ in $I\Delta_0+\Omega_1$ by finite rooted acyclic graphs;
bisimulation supplies equality, and bounded truth on the graphs supplies
$\Delta_0$-Separation.  We then characterize fragments of set-theoretic
Collection by Gaifman-style splitting and cofinal elementarity.  Finally,
we separate two end-extension mechanisms.  Without a resolution,
Kaufmann's construction characterizes the relevant Collection fragments
for countable and locally for $\aleph_1$-like models.  With a strict
transitive resolution, Power Set, Infinity, and strong Collection, a
finite-strength Keisler--Morley construction gives $\Sigma_N$-elementary
taller* end extensions whose output retains strong
$\Sigma_{N-2}$-Collection.  The proof is choice-free inside the model and
also recovers the classical $\mathsf{ZF}$ theorem.
\end{abstract}

\keywords{collection; end extension; cofinal extension; Gaifman splitting;
interpretability; weak arithmetic; weak set theory}

\section{Introduction}

Collection has parallel formulations in arithmetic and set theory.  For a
formula $\varphi$, an arithmetic instance has the form
\[
 \forall v\forall p\bigl(\forall x<p\,\exists y\,\varphi(x,y,v)
 \longrightarrow \exists q\forall x<p\,\exists y<q\,\varphi(x,y,v)\bigr),
\]
whereas the corresponding set-theoretic instance replaces the initial
segment $p$ by the members of a set $p$:
\[
 \begin{aligned}
 \forall v\forall p\bigl(&\forall x\in p\,\exists y\,\varphi(x,y,v)\\
 &\longrightarrow \exists q\forall x\in p\,\exists y\in q\,\varphi(x,y,v)\bigr).
 \end{aligned}
\]
This formal analogy raises a structural question: how much of the model
theory of weak arithmetic survives when bounded initial segments are
replaced by sets?  The answer depends sharply on the base theory.  In
arithmetic, $\mathsf{PA}^-$ is an algebraic theory, substantially weaker
than the induction theories $I\Sigma_n$, but it already supports robust
coding and end-extension arguments.  Mathias's theory $\mathsf{DB}_0$
plays a comparable role on the set-theoretic side: it has the elementary
rudimentary operations and $\Delta_0$-Separation, but no Infinity, Power
Set, Foundation, or Collection \cite{Ma06}.

Gaifman's splitting theorem says that if
$\mathcal{M}\preccurlyeq\mathcal{N}$ are models of $\mathsf{PA}$, then the
substructure whose domain is
$\operatorname{sup}_{\mathcal N}(M)=\{b\in N:\exists a\in M\ (b<a)\}$
lies elementarily between them \cite{Ga72}.  Kaye, and more recently
Kurahashi and Minami, used variants of this property to characterize
induction and collection schemes over weak arithmetic
\cite{Ka91,Ku24}.  The set-theoretic analogue replaces the supremum by the
membership hull
\[
 M^*=\{b\in N:\exists a\in M\ (\mathcal N\vDash b\in a)\}.
\]
We develop corresponding characterizations over $\mathsf{DB}_0$.

Strong
$\Sigma_{n+1}$-Collection is exactly the assertion that the membership
hull splits every suitable extension at level $\Sigma_{n+1}$.  Ordinary
$\Sigma_{n+1}$-Collection is recovered by combining a lower splitting
condition with cofinal elementarity.  At all finite levels these results
yield the full Gaifman splitting property precisely when full Collection
holds.

Before turning to splitting, we make the comparison with weak arithmetic
precise in two independent ways.  First, every $\mathsf{PA}^-$-model and
every $\mathsf{DB}_0$-model has a proper end extension satisfying the same
base theory.  The set-theoretic construction is external and in fact adds
an exact transitive cover of the old universe.  Second, and more
substantially, we prove
\[
   \mathsf{DB}_0\equiv_{\mathrm{int}}\mathsf{PA}^-.
\]
The difficult direction factors through $I\Delta_0+\Omega_1$.  Its
objects are finite rooted directed acyclic graphs with explicit padding
bounds.  Equality is bisimulation rather than equality of graph codes.
This choice is essential: ordinary Ackermann coding would require the
uniform existence of $2^a$, while $\Omega_1$ controls only fixed
polynomials in the binary length of the input.  Local satisfaction on the
graphs then makes every translated $\Delta_0$-Separation instance
arithmetically bounded.

The final part concerns end extensions.  The MacDowell--Specker theorem
gives every model of $\mathsf{PA}$ an elementary end extension \cite{Do90}.
Keisler and Morley proved the set-theoretic analogue for models of
$\mathsf{ZF}$ whose ordinals have countable cofinality \cite{Ke68}, and
Kaufmann related partial elementary end extensions of weak set theories to
fragments of Collection \cite{Ka81}.  Two mechanisms should be kept
separate.  Kaufmann's omitting-types argument needs no resolution and gives
the countable, as well as an $\aleph_1$-like local, characterization.  The
Keisler--Morley construction uses rank stages as a resolution.  We replace
them by a definable strict transitive resolution and determine one explicit
finite base theory sufficient for the construction.  If $N\geq6$, the
input theory is
\[
 \B_N=\mathsf{DB}_0+\mathsf{PowerSet}+\mathsf{Infinity}
      +\mathsf{Coll}_s(\Sigma_N),
\]
and the resulting $\Sigma_N$-elementary model satisfies $\B_{N-2}$.  The
proof uses set-valued Skolem envelopes and the choice-free finite
Erd\H{o}s--Rado argument from Keisler and Morley's appendix; it does not
assume Choice inside the model.

Section~2 fixes all common syntactic, interpretability, extension, and
resolution terminology.  Section~3 compares $\mathsf{DB}_0$ with
$\mathsf{PA}^-$ by end extensions and interpretations.  Section~4 proves
the Gaifman-style characterizations.  Section~5 first treats end extensions
without resolutions and then proves the finite-strength Keisler--Morley
theorem under a strict transitive resolution.  Section~6 records the
remaining questions.

\section{Preliminaries}

We work in the language $\mathcal L_\in=\{\in\}$, with equality treated as
logical.  A quantifier is \emph{bounded} if it has the form
$\exists x\in y$ or $\forall x\in y$.  The classes $\Delta_0=\Sigma_0=
\Pi_0$ and $\Sigma_n,\Pi_n$ are the usual Levy hierarchy; finite tuples of
parameters are allowed.  Formula classes are understood up to provable
equivalence over the displayed base theory.

If $\mathcal M\subseteq\mathcal N$, then
$\mathcal M\preccurlyeq_\Gamma\mathcal N$ means that the two structures
agree on all $\Gamma$-formulas with parameters from $M$.  We omit the
subscript for full elementarity.  In the language
$\mathcal L_\in(M)=\mathcal L_\in\cup\{c_a:a\in M\}$, the elementary
diagram $\ElDiag(\mathcal M)$ is the set of all sentences true in
$(\mathcal M,a)_{a\in M}$.  We similarly write
$\operatorname{Diag}_\Gamma(\mathcal M)$ for its $\Gamma$-part.  These
conventions also fix the meaning of diagram notation used in the compactness
arguments below.

For theories $T$ and $U$, we write $T\triangleright U$ when $T$ interprets
$U$ in the usual first-order sense: the interpretation may use a definable
domain and a definable congruence for equality, followed by passage to the
quotient.  Fixed finite-dimensional interpretations are identified with
one-dimensional interpretations by a standard tuple code.  We write
$T\equiv_{\mathrm{int}}U$ for mutual interpretability; no definable
isomorphism between the two composites and the identity interpretations is
implicit in this notation.

We use the usual languages for weak arithmetic, up to harmless definitional
extensions.  The theory $\mathsf Q$ is Robinson arithmetic, and
$\mathsf{PA}^-$ is the theory of the nonnegative parts of discretely ordered
commutative semirings.  The theory $I\Delta_0$ has induction for bounded
formulas.  If $|x|$ denotes binary length, our version of $\Omega_1$ asserts
the totality of
\[
   \omega_1(x)=2^{|x|^2}.
\]
We use the bounded bit predicate $\Bit(a,i)$ and the standard partial
relation $\PowTwo(p,k)$ for $p=2^k$.  These are conservative definitional
devices: $\PowTwo(p,k)$ can be expressed by a bounded doubling table once
an $\Omega_1$ bound is supplied.  In particular, we never assume that the
full exponential function is total.

Following \cite{Ma06,Mc19,Ku24}, let $\mathsf{DB}_0$ consist of
Extensionality, Empty Set, Pairing, Union, Cartesian Product, and the
$\Delta_0$-Separation scheme.  For a formula
$\varphi(x,y,\bar v)$, put
\begin{align*}
 \mathsf{Coll}_s(\varphi):\quad
 &\forall\bar v\,\forall p\,\exists q\,\forall x\in p
   \bigl(\exists y\,\varphi(x,y,\bar v)
   \leftrightarrow
   \exists y\in q\,\varphi(x,y,\bar v)\bigr),\\
 \mathsf{Coll}(\varphi):\quad
 &\forall\bar v\,\forall p
   \bigl(\forall x\in p\,\exists y\,\varphi(x,y,\bar v)
   \longrightarrow
   \exists q\,\forall x\in p\,\exists y\in q\,
   \varphi(x,y,\bar v)\bigr),\\
 \mathsf{Coll}_w(\varphi):\quad
 &\forall\bar v
   \bigl(\forall x\,\exists y\,\varphi(x,y,\bar v)
   \longrightarrow
   \forall p\,\exists q\,\forall x\in p\,\exists y\in q\,
   \varphi(x,y,\bar v)\bigr).
\end{align*}
For $\theta(x,\bar v)$, let
\[
 \mathsf{Sep}(\theta):\quad
 \forall\bar v\,\forall p\,\exists q\,\forall x
 \bigl(x\in q\leftrightarrow
       (x\in p\wedge\theta(x,\bar v))\bigr).
\]
A superscript ${}^{-}$ denotes the parameter-free version: delete
$\bar v$ and its universal quantifier.  Thus, for example,
$\mathsf{Coll}^{-}(\varphi)$ applies to a formula $\varphi(x,y)$, not to
a formula with an unquantified parameter.  For a syntactic class $\Gamma$,
$\mathsf{Coll}(\Gamma)$ denotes the scheme containing the instances for
all $\Gamma$-formulas, and similarly for the other schemes.

We first record the elementary relations between these variants.

\begin{theorem}[McKenzie,\cite{Mc19}]\label{thm:strong-collection-equivalence}
    For $n\geq0$, the theory
    $\mathsf{DB}_0+\mathsf{Coll}_s(\Sigma_{n+1})$ is deductively
    equivalent to
    \[
       \mathsf{DB}_0+\mathsf{Coll}(\Sigma_{n+1})
       +\mathsf{Sep}(\Sigma_{n+1}).
    \]
\end{theorem}
\begin{proof}
    McKenzie's proof uses only the axioms listed in the definition of
    $\mathsf{DB}_0$.  In particular, Infinity, Set Foundation, Power Set,
    and transitive containment play no role, so the proof applies verbatim
    to the present base theory.
\end{proof}

\begin{theorem}\label{thm:collection-variants}
    Let $n\geq0$.  Over $\mathsf{DB}_0$, the eight schemes obtained by
    choosing one of
    \[
       \mathsf{Coll},\quad \mathsf{Coll}^{-},\quad
       \mathsf{Coll}_w,\quad \mathsf{Coll}_w^{-}
    \]
    and one of the formula classes $\Sigma_{n+1}$ and $\Pi_n$ are
    pairwise equivalent.  Moreover,
    $\mathsf{Coll}_s(\Sigma_{n+1})$ and
    $\mathsf{Coll}_s(\Pi_n)$ are equivalent.
\end{theorem}
\begin{proof}
    We indicate the three standard reductions.  First, existential witnesses
    in a $\Sigma_{n+1}$ formula may be included in the collected object,
    while a $\Pi_n$ formula may be regarded as a $\Sigma_{n+1}$ formula by
    adding a dummy existential quantifier.  Pairing and two applications of
    Union recover the original witness from its code.  This proves the
    equivalence between the two complexity classes.

    Second, to apply weak Collection to a set $p$, replace
    $\varphi(x,y,\bar v)$ by
    \[
       x\notin p\ \vee\ \varphi(x,y,\bar v).
    \]
    Under the antecedent of ordinary Collection this relation is total, and
    the displayed modification does not raise Levy complexity.  The converse
    implication is immediate.

    Finally, parameters are coded into the first argument.  Apply a
    parameter-free instance to $p\times\{\bar v\}$ and to the formula which,
    on an ordered pair $\langle x,\bar v\rangle$, expresses
    $\varphi(x,y,\bar v)$.  The standard pairing operations are
    $\mathsf{DB}_0$-suitable \cite[Definition~1.4 and
    Proposition~2.64]{Ma06}, so this coding preserves Levy complexity.
    The same reductions work for strong Collection.
\end{proof}

The following bounded-quantifier calculation will be used repeatedly.  If
$\Gamma$ is a formula class, write $\Delta_0(\Gamma)$ for formulas obtained
from $\Gamma$-formulas by finite conjunctions, disjunctions, and bounded
quantifiers, without taking complements of the distinguished
$\Gamma$-subformulas.

\begin{lemma}[Bounded contraction]
\label{lem:bounded-contraction}
Let $n\geq1$.  Over
$\mathsf{DB}_0+\mathsf{Coll}(\Sigma_n)$, every
$\Delta_0(\Sigma_n)$ formula is equivalent to a $\Sigma_n$ formula and
every $\Delta_0(\Pi_n)$ formula is equivalent to a $\Pi_n$ formula.  The
equivalences are uniform in all displayed parameters.
\end{lemma}
\begin{proof}
Proceed by induction on the bounded syntax.  The only clauses not already
absorbed by the usual closure properties are a bounded universal quantifier
over a $\Sigma_n$ formula and, dually, a bounded existential quantifier over
a $\Pi_n$ formula.  Write the first formula as
\[
   \forall x\in a\,\exists z\,\theta(x,z,\bar v),
   \qquad \theta\in\Pi_{n-1}.
\]
By Theorem~\ref{thm:collection-variants},
$\mathsf{Coll}(\Sigma_n)$ includes Collection for $\Pi_{n-1}$ formulas, so
the displayed formula is equivalent to
\[
   \exists b\,\forall x\in a\,\exists z\in b\,
      \theta(x,z,\bar v).                                  \tag{C}
\]
The matrix after $\exists b$ is in $\Pi_{n-1}$ after applying the induction
hypothesis to its bounded quantifiers.  Thus (C) is $\Sigma_n$.  Negating
this equivalence contracts a bounded existential over a $\Pi_n$ formula to
$\Pi_n$.  Conjunctions, disjunctions, and the remaining bounded quantifiers
are handled inductively.  No Collection is needed for the corresponding
$n=0$ calculation, since all formulas involved are already bounded.
\end{proof}

\begin{theorem}[Schindler--Schlicht]\label{thm:parameter-separation}
    For $n\geq0$, the schemes $\mathsf{Sep}(\Sigma_{n+1})$ and
    $\mathsf{Sep}^{-}(\Sigma_{n+1})$ are equivalent over
    $\mathsf{DB}_0$.
\end{theorem}
\begin{proof}
    The parameter-coding proof of \cite[Lemma~0.4]{Sc16} requires only the
    sets of one- and two-element subsets of a given set, rather than its full
    power set.  Their existence is provable in $\mathsf{DB}_0$
    \cite[Proposition~2.60]{Ma06}.
\end{proof}

\begin{theorem}\label{thm:weak-to-strong}
    If $n\geq1$, then
    \[
       \mathsf{DB}_0+\mathsf{Coll}_w^{-}(\Sigma_{n+1})
       \vdash\mathsf{Coll}_s(\Sigma_n).
    \]
    For $n=0$, the same argument gives
    $\mathsf{Coll}_w^{-}(\Sigma_2)\vdash
      \mathsf{Coll}_s(\Delta_0)$.
\end{theorem}
\begin{proof}
    Let $\varphi(x,y,\bar v)$ be $\Sigma_n$, and fix $p,\bar v$.
    The relation
    \[
       \psi(\langle x,\bar v\rangle,z):=
       \bigl(\exists y\,\varphi(x,y,\bar v)
       \longrightarrow\varphi(x,z,\bar v)\bigr)
    \]
    is total.  For $n\geq1$ it is equivalent to a
    $\Sigma_{n+1}$ formula, because Boolean combinations of
    $\Sigma_n$ and $\Pi_n$ formulas lie at the next level.  Apply weak
    parameter-free Collection to $p\times\{\bar v\}$.  The resulting set
    contains a $\varphi$-witness for every $x\in p$ for which one exists,
    which is precisely the nontrivial direction of strong Collection.
    When $n=0$, the unbounded existential quantifier in the antecedent makes
    $\psi$ a $\Pi_1$ formula, hence a $\Sigma_2$ formula; this explains the
    separately stated base case.
\end{proof}

\begin{theorem}\label{thm:strong-parameters}
    For every $n\geq0$, the schemes
    $\mathsf{Coll}_s^{-}(\Sigma_n)$ and
    $\mathsf{Coll}_s(\Sigma_n)$ are equivalent over $\mathsf{DB}_0$.
\end{theorem}
\begin{proof}
    Only parameter elimination needs proof.  Given
    $\varphi(x,y,\bar v)$ and $p,\bar v$, apply the parameter-free instance
    to $p\times\{\bar v\}$ and to the formula obtained by decoding the first
    argument as $\langle x,\bar v\rangle$.  Suitability of pairing preserves
    the complexity, as in Theorem~\ref{thm:collection-variants}.
\end{proof}

We now fix the extension terminology used throughout the paper.  For
$a\in M$ let
\[
   \operatorname{Ext}_{\mathcal M}(a)
   :=\{x\in M:\mathcal M\vDash x\in a\}.
\]
If $\mathcal M\subseteq\mathcal N$, then $\mathcal N$ is an \emph{end
extension} of $\mathcal M$, written
$\mathcal M\subseteq_{\mathsf{end}}\mathcal N$, if
\[
   \operatorname{Ext}_{\mathcal N}(a)
   =\operatorname{Ext}_{\mathcal M}(a)
   \qquad(a\in M).
\]
It is a \emph{cofinal extension}, written
$\mathcal M\subseteq_{\mathsf{cf}}\mathcal N$, if every $b\in N$ is a
member, in $\mathcal N$, of some $a\in M$.  An element $a\in M$ is
\emph{fixed} by $\mathcal N$ when the displayed equality of extensions
holds for that particular $a$.

The extension is \emph{taller} if there is $c\in N$ such that
$\mathcal N\vDash a\in c$ for every $a\in M$.  It is \emph{taller*} if
such a $c$ can be chosen transitive in $\mathcal N$.  We call $c$ an
\emph{exact cover} of $M$ if
$\operatorname{Ext}_{\mathcal N}(c)=M$.

Given $\mathcal M\preccurlyeq\mathcal N$, define its membership hull by
\[
  M^*=\{b\in N:\exists a\in M\ (\mathcal N\vDash b\in a)\},
  \qquad
  \Mstar=(M^*,\in^{\mathcal N}\mathbin{\upharpoonright}(M^*)^2).
\]
Pairing gives $M\subseteq M^*$.  By definition
$\mathcal M\subseteq_{\mathsf{cf}}\Mstar$; using Union once shows that
$\Mstar\subseteq_{\mathsf{end}}\mathcal N$.

The axiom of transitive containment, denoted by $\mathsf{TCo}$, is
\[
  \forall x\,\exists t\bigl(x\in t\wedge\Trans(t)\bigr),
  \qquad
  \Trans(t):=\forall u\in t\,\forall v\in u\,(v\in t).
\]
It is a $\Pi_2$ sentence.  It is weaker than $\mathsf{TCl}$, the assertion
that every set has a least transitive superset.  The finite-chain argument
relating the two principles, including the precise role of Set Foundation,
is recorded in Remark~\ref{rem:tco-and-tcl}; none of the main extension
theorems uses $\mathsf{TCl}$.

We next isolate the substitute for a rank filtration used in the last
section.  A \emph{resolution of complexity $r$} of $\mathcal M$ consists of
a parameter $v\in M$ and formulas
\[
  \rho_\Sigma(P,\alpha,v)\in\Sigma_r,
  \qquad
  \rho_\Pi(P,\alpha,v)\in\Pi_r,
\]
which are equivalent in $\mathcal M$ and uniquely define a stage
$P_\alpha$ for every $\alpha\in\Ord^{\mathcal M}$.  The stages are
increasing and cover the universe:
\[
 \alpha\in\beta\Longrightarrow P_\alpha\subseteq P_\beta,
 \qquad
 \forall x\,\exists\alpha\in\Ord\ (x\in P_\alpha).
\]
Thus stage membership
\[
 R(x,\alpha):\Longleftrightarrow
 \exists P\,(\rho_\Sigma(P,\alpha,v)\wedge x\in P)
\]
has both a $\Sigma_r$ and a $\Pi_r$ definition.  A resolution is
\emph{strict transitive} if, in addition, its ordinal indices are linearly
ordered by membership and contain zero, successors, and finite suprema,
each $P_\alpha$ is transitive and $P_\alpha\notin P_\alpha$, and
\[
    P_\alpha\in P_\beta
       \quad\Longleftrightarrow\quad \alpha\in\beta.       \tag{R}
\]
Condition (R) is what permits the stages themselves to realize arbitrary
finite linear orders.  The rank hierarchy in a model of $\mathsf{ZF}$ is
the motivating example.

For $n\geq1$, set
\[
  \B_n:=\mathsf{DB}_0+\mathsf{PowerSet}+\mathsf{Infinity}
          +\mathsf{Coll}_s(\Sigma_n).
\]
Power Set supplies spaces of colors, Infinity supplies arbitrarily large
finite partition bounds, and strong Collection supplies the corresponding
Separation and functional Replacement.

Finally, a structure $\mathcal M$ is \emph{$\aleph_1$-like} if
$|M|=\aleph_1$ externally while
$\operatorname{Ext}_{\mathcal M}(a)$ is countable for every $a\in M$.

\begin{definition}
    Let $\mathcal M\vDash\mathsf{DB}_0$ and $n\geq0$.
    \begin{enumerate}
        \item $\mathcal M$ has $\mathbf{end}_n$ if
        $\Mstar\preccurlyeq_{\Sigma_n}\mathcal N$ for every
        $\mathcal M\preccurlyeq\mathcal N$.
        \item $\mathcal M$ has $\mathbf{cof}_n$ if
        $\mathcal M\preccurlyeq_{\Sigma_n}\Mstar$ for every
        $\mathcal M\preccurlyeq\mathcal N$.
        \item $\mathcal M$ has $\mathbf{COF}_n$ if, whenever
        $\mathcal M\preccurlyeq_{\Delta_0}\mathcal K$,
        $\mathcal M\subseteq_{\mathsf{cf}}\mathcal K$, and
        $\mathcal K\vDash\mathsf{DB}_0$, one has
        $\mathcal M\preccurlyeq_{\Sigma_n}\mathcal K$.
    \end{enumerate}
\end{definition}

\section{Comparing \texorpdfstring{$\mathsf{DB}_0$}{DB0} and
\texorpdfstring{$\mathsf{PA}^-$}{PA-minus}}

\subsection{Proper end extensions}

We begin with a direct structural parallel.  It will later be complemented
by a comparison at the level of interpretations.

\begin{theorem}[{\cite[Exercise~7.7]{Ka91a}}]
\label{thm:pa-minus-end-extension}
    Every model of $\mathsf{PA}^-$ has a proper end extension which is again
    a model of $\mathsf{PA}^-$.
\end{theorem}

\begin{theorem}\label{thm:db0-exact-cover}
    Every model of $\mathsf{DB}_0$ has a proper taller* end extension
    which is again a model of $\mathsf{DB}_0$.  Moreover, the taller*
    witness can be chosen to be an exact cover of the original model.
\end{theorem}
\begin{proof}
    We use an external closure construction.  If
    $\mathcal{A}=(A,E)$ is any extensional structure, let
    $\operatorname{Def}(\mathcal{A})$ be the family of subsets of $A$
    that are first-order definable in $\mathcal{A}$ with parameters.  Form
    an extensional structure $\mathcal{A}^{+}$ as follows.  For each
    $S\in\operatorname{Def}(\mathcal{A})$ which is not of the form
    \[
        \operatorname{Ext}_{\mathcal{A}}(a)
        =\{x\in A: \mathcal{A}\vDash x\mathrel{E}a\}
        \qquad(a\in A),
    \]
    add one new element $s_S$.  Retain $E$ on $A$, put
    $x\mathrel{E^{+}}s_S$ exactly when $x\in S$, and put no new element
    into the extension of an old element.  Distinct new elements have
    distinct extensions, and none has the same extension as an old
    element.  Hence $\mathcal{A}^{+}$ is extensional and
    $\mathcal{A}\subseteq_{\mathsf{end}}\mathcal{A}^{+}$.

    Starting with $\mathcal{M}_0=\mathcal{M}$, recursively set
    $\mathcal{M}_{k+1}=\mathcal{M}_k^{+}$ and let
    \[
        \mathcal{N}=\bigcup_{k<\omega}\mathcal{M}_k.
    \]
    Each map in this chain is inclusion and is an end embedding.  It
    follows that $\mathcal{M}\subseteq_{\mathsf{end}}\mathcal{N}$ and
    that $\mathcal{N}$ is extensional.

    We verify the remaining axioms of $\mathsf{DB}_0$.  The empty set is
    already in $M$.  If $a,b\in N$, choose $k$ with $a,b\in M_k$.  The
    subset $\{a,b\}$ of $M_k$ is definable over $\mathcal{M}_k$, so it is
    represented in $M_{k+1}$.  Similarly, if $a\in M_k$, then
    \[
        \{x\in M_k:\mathcal{M}_k\vDash
          \exists y\in a\,(x\in y)\}
    \]
    is represented at the next stage and is the union of $a$ in
    $\mathcal{N}$.  This remains the correct union because later stages
    add no new members to $a$ or to any member of $a$.  Repeated use of
    pairing produces, uniformly within three further stages, every
    Kuratowski pair $\langle x,y\rangle$ with $x\in a$ and $y\in b$.
    The set of all these pairs is then a definable subset of that stage,
    so one further application of the construction produces $a\times b$.

    Finally, let $a,\bar p\in M_k$ and let $\theta(x,\bar p)$ be
    $\Delta_0$.  Bounded formulas are absolute between an end
    substructure and its union, so
    \[
        S=\{x\in M_k:\mathcal{M}_k\vDash
        x\in a\wedge\theta(x,\bar p)\}
         =\{x\in N:\mathcal{N}\vDash
        x\in a\wedge\theta(x,\bar p)\}.
    \]
    The set $S$ is represented in $M_{k+1}$.  Thus $\mathcal{N}$
    satisfies $\Delta_0$-separation, and hence all the axioms of
    $\mathsf{DB}_0$ in the formulation of \cite{Ma06}.

    The full domain $M$, viewed externally as the subset defined by
    $x=x$, is not the extension of any old element.  For if some $a\in M$
    contained every element of
    $M$, $\Delta_0$-separation on $a$ would produce
    $r=\{x\in a:x\notin x\}$.  Since $r\in M$, universality of $a$ would
    give $r\in a$, and hence
    $r\in r\leftrightarrow r\notin r$, a contradiction.  Thus the first
    stage adds an element $c=s_M$ whose members are exactly the elements
    of $M$.  Later stages add no members to $c$, so
    $\operatorname{Ext}_{\mathcal N}(c)=M$ and $c\notin c$.  Moreover, $c$
    remains transitive: if $y\in x\in c$, then $x\in M$ and endness puts
    $y$ in $M$, hence in $c$.  Therefore $c$ witnesses that the proper end
    extension $\mathcal{N}$ is taller*.
\end{proof}

\begin{remark}
    The construction is external and proves only existence.  In particular,
    it does not produce the kind of canonical minimal end extension available
    for $\mathsf{PA}^-$.  Non-well-founded end extensions may disagree even
    on whether a new cover belongs to itself, so any set-theoretic minimality
    statement would require additional hypotheses and a specified category of
    embeddings.
\end{remark}

\subsection{Mutual interpretability}

Let $\mathsf{AS}$ denote extensional adjunctive set theory: Extensionality,
Empty Set, and the axiom saying that $x\cup\{y\}$ exists for every $x,y$.
\begin{proposition}\label{prop:db0-interprets-q}
    The theory $\mathsf{DB}_0$ interprets Robinson arithmetic
    $\mathsf Q$.
\end{proposition}
\begin{proof}
    Empty Set, Pairing, and Union show that $\mathsf{DB}_0$ extends
    $\mathsf{AS}$.  The Szmielew--Tarski interpretation, in the form
    proved in \cite{Co70} and presented in
    \cite[Sections~3.1--3.2]{Vis08}, interprets
    $\mathsf Q$ in $\mathsf{AS}$.  Composing with the identity
    interpretation of $\mathsf{AS}$ in $\mathsf{DB}_0$ gives the result.
\end{proof}
\begin{remark}
    The converse direction requires more than the known interpretation of
    adjunctive set theory in $\mathsf Q$.
    Union, Cartesian Product, and the whole $\Delta_0$-Separation scheme
    still have to be verified.  The graph interpretation below supplies
    precisely these missing operations.
\end{remark}

Put
\[
   S:=I\Delta_0+\Omega_1.
\]
Wilkie's definable-cut construction gives
$\mathsf Q\triangleright S$ \cite{HP93}; see also
\cite[p.~307]{Vis08}.  We now construct the missing interpretation
$S\triangleright\mathsf{DB}_0$.  Finite acyclic
graph representations of hereditarily finite sets are classical; compare
\cite{Da17,Sa93}.  The explicit padding and bounded-truth argument below
is needed to formalize all of $\mathsf{DB}_0$ in the weak theory $S$.

\subsubsection{Short computations in \texorpdfstring{$S$}{S}}

We work in the standard conservative coding expansion of $S$ containing
binary length, bits, fixed tuple codes, and finite-table codes.  If these
symbols are eliminated, the finitely many $\Omega_1$ witnesses used in a
given argument are first quantified and then used as bounds.  All induction
formulas below remain bounded formulas in the original language.

\begin{lemma}[Short computations]\label{lem:short-computations}
Let
\[
    L=1+|x_0|+\cdots+|x_{r-1}|
\]
be the total binary length of a fixed-length tuple.
\begin{enumerate}
    \item For every standard $d,c<\omega$, the theory $S$ proves that
    there are $K,P$ such that
    \[
       \PowTwo(P,K)\quad\text{and}\quad K\geq cL^d.       \tag{I1}
    \]
    \item If $\PowTwo(P,N)$ and $\theta(i,\bar a)$ is a fixed bounded
    arithmetic formula, then there is $b<P$ such that
    \[
       \forall i<N\,
       \bigl(\Bit(b,i)\leftrightarrow\theta(i,\bar a)\bigr). \tag{I2}
    \]
    The same holds for a bit string whose $i$th bit is uniquely determined
    by its earlier bits through a fixed bounded formula.
    \item Every fixed clocked deterministic binary algorithm whose time,
    space, and output length are bounded by a fixed polynomial in $L$ has,
    provably in $S$, a unique computation table and output on every input.
\end{enumerate}
\end{lemma}
\begin{proof}
For (1), let $b_0$ be a fixed tuple code for the inputs, enlarged by a
constant, and iterate $b_{i+1}=\omega_1(b_i)$ a standard finite number of
times.  Up to a fixed constant factor, $|b_0|$ is $L$, while
\[
   |b_{i+1}|\geq |b_i|^2.
\]
Choose a standard $e$ with $2^e>d$, increasing it to absorb $c$ and the
finitely many small inputs.  Taking $K=|b_{e-1}|^2$ and $P=b_e$ gives
(I1).  The number of iterations depends only on $c,d$, not on a number
quantified inside $S$.

For (2), bounded induction on $i\leq N$ first supplies all initial powers
$2^i$, since $2^N$ exists.  Starting with $b_0=0$, put
$b_{i+1}=b_i+2^i$ when $\theta(i,\bar a)$ holds and put
$b_{i+1}=b_i$ otherwise.  Bounded induction constructs the code and proves
its correctness.  A recursion which reads only earlier bits is handled by
the same prefix construction.

For (3), suppose time and space are at most $L^d$.  A complete computation
table has at most $CL^{2d}$ bits for a standard $C$.  Part (1) supplies a
power of two large enough to code it.  The local transition rule is
bounded, so bounded prefix induction constructs the normalized table and
proves uniqueness.  The output is then extracted by a bounded formula.
\end{proof}

The point of the lemma is local.  Each fixed translated
$\Delta_0$-Separation instance may use its own polynomial and its own fixed
number of $\omega_1$ iterations; interpreting a scheme does not require a
uniform truth predicate or total exponentiation.

\subsubsection{Rooted DAGs and bisimulation}

An object of the interpretation is coded by a tuple
\[
   G=(n,k,p,A,r).
\]
Fixed arithmetic tuple coding makes this a one-dimensional interpretation.
Let $D(G)$ say that
\[
\begin{gathered}
  1\leq n,\qquad r<n,\qquad n^2\leq k,\qquad
  \PowTwo(p,k),\qquad A<p,\\
  \forall i,j<n\,
  \bigl(\operatorname{Edge}_G(i,j)\longrightarrow j<i\bigr),
\end{gathered}\tag{I3}
\]
where
\[
   \operatorname{Edge}_G(i,j)
       :\Longleftrightarrow \Bit(A,in+j).
\]
Only the first $n^2$ bits of $A$ are read.  The field $p=2^k$ is an
explicit padding bound: even an edgeless graph carries enough room to code
relations on all of its vertices.  Edges point to smaller indices, hence
the graph is acyclic.  For $i<n$, write
\[
   G[i]=(n,k,p,A,i)
\]
for the same graph re-rooted at $i$.  Vertices not reachable from the root
are harmless padding.

Let $G=(n,k,p,A,r)$ and $H=(m,\ell,q,B,s)$ satisfy $D$.  A bisimulation
table has entries $C(i,j)$, for $i<n$ and $j<m$, determined by
\[
\begin{aligned}
C(i,j)\ \longleftrightarrow\ {}
&\forall i'<n\bigl(\operatorname{Edge}_G(i,i')\rightarrow
   \exists j'<m\,(\operatorname{Edge}_H(j,j')\wedge C(i',j'))\bigr)\\
{}\wedge{}&\forall j'<m\bigl(\operatorname{Edge}_H(j,j')\rightarrow
   \exists i'<n\,(\operatorname{Edge}_G(i,i')\wedge C(i',j'))\bigr).
\end{aligned}\tag{I4}
\]
The right side reads only entries with smaller indices.  Moreover,
\[
   nm\leq\max(n^2,m^2)\leq\max(k,\ell),
\]
so the input padding bounds code the entire table.  By
Lemma~\ref{lem:short-computations}(2), $S$ proves that there is a unique
normalized table satisfying (I4).  Define
\[
  G[i]\sim H[j]\quad:\Longleftrightarrow\quad C(i,j),
  \qquad
  G\approx H\quad:\Longleftrightarrow\quad G[r]\sim H[s]. \tag{I5}
\]

\begin{lemma}\label{lem:bisimulation-congruence}
The theory $S$ proves that $\approx$ is an equivalence relation on $D$.
If
\[
 G\mathrel EH\quad:\Longleftrightarrow\quad
 \exists j<m\bigl(\operatorname{Edge}_H(s,j)\wedge G\approx H[j]\bigr),
                                                               \tag{I6}
\]
then $E$ is invariant under $\approx$ in both coordinates.
\end{lemma}
\begin{proof}
Reflexivity follows by bounded induction on vertex indices, and symmetry
interchanges the two halves of (I4).  For transitivity, suppose
$G[i]\sim H[j]$ and $H[j]\sim K[t]$.  Bounded induction on $i+j+t$
composes the two back-and-forth matches: every matching successor has
strictly smaller indices, where the induction hypothesis applies.  Thus
$G[i]\sim K[t]$.  For invariance of (I6), use the back-and-forth clause at
the root of $H\approx H'$ to match a root successor of $H$ with one of
$H'$, and then use symmetry and transitivity to replace $G$ by an
equivalent $G'$.
\end{proof}

Consequently $E$ is a well-defined relation on $D/{\approx}$.  Arithmetic
equality of graph codes would not work here: duplicated successors,
sharing, and unreachable padding must all represent the same set.

\begin{lemma}[Local bounded truth]\label{lem:local-bounded-truth}
For every set-theoretic $\Delta_0$ formula $\varphi(\bar x)$ there is a
bounded arithmetic formula $\operatorname{Sat}_\varphi(\bar G)$ such that
$S$ proves
\[
   \varphi^I(\bar G)
      \longleftrightarrow\operatorname{Sat}_\varphi(\bar G), \tag{I7}
\]
where $I=(D,\approx,E)$ is the interpretation just defined.  This bounded
satisfaction relation is invariant when any $G_i$ is replaced by an
$\approx$-equivalent code.
\end{lemma}
\begin{proof}
Define satisfaction externally by induction on $\varphi$.  For atomic
formulas use $\approx$ and $E$, and translate Boolean connectives
literally.  A bounded existential is localized to the root successors:
\[
\begin{split}
 \operatorname{Sat}_{\exists x\in y\,\psi}(H,\bar G)
 \quad:\Longleftrightarrow\quad
 \exists j<m_H\bigl(&\operatorname{Edge}_H(r_H,j)\\
   &{}\wedge\operatorname{Sat}_\psi(H[j],\bar G)\bigr).
\end{split}\tag{I8}
\]
All quantifiers in (I8) are bounded by the finite input graph; the
bisimulation tables occurring in the atomic cases are bounded by the input
padding.  Formula induction proves (I7).  In its quantified step, the
back-and-forth clauses match the root successors of equivalent bound
objects.  The same induction therefore proves invariance.
\end{proof}

We shall use one construction principle repeatedly.  Suppose that a new
graph has a number of vertices, and hence a number of adjacency bits,
bounded by a fixed polynomial in the binary lengths of the input graph
codes, and each new edge is decided by a fixed bounded formula or short
computation.  Lemma~\ref{lem:short-computations} supplies $K,P$ with
$P=2^K$ and $K$ above the square of the new vertex count; bounded bit
comprehension then supplies its adjacency matrix.  Thus $S$ proves that the
new padded DAG code exists.  We refer to this as the graph-construction
principle.

\begin{theorem}\label{thm:dag-interpretation}
The theory $I\Delta_0+\Omega_1$ interprets $\mathsf{DB}_0$.
\end{theorem}
\begin{proof}
We verify the axioms in the quotient $D/{\approx}$.

\emph{Extensionality.}
Suppose $G$ and $H$ have the same $E$-members.  For every root successor
$i$ of $G$, the re-rooted code $G[i]$ is a member of $G$, hence a member of
$H$.  It is therefore bisimilar to a root successor of $H$.  The converse
is symmetric, so (I4) gives $G\approx H$.

\emph{Empty Set and Pairing.}
The code $(1,1,2,0,0)$ has no members.  Given $G,H$, take disjoint copies
and add a new largest root pointing to the two old roots.  Call the result
$P(G,H)$.  The graph-construction principle gives its code, and
\[
   X\mathrel E P(G,H)
      \quad\longleftrightarrow\quad
      (X\approx G\ \vee\ X\approx H).                    \tag{I9}
\]
Thus the construction also works when $G\approx H$.

\emph{Union.}
Copy the graph of $H$, add a new root, and point it to $j$ exactly when
there is a root successor $i$ of $H$ with an edge from $i$ to $j$.  If the
result is $U(H)$, then invariance from
Lemma~\ref{lem:bisimulation-congruence} gives
\[
 X\mathrel E U(H)
 \quad\longleftrightarrow\quad
 \exists Y\in D\,(Y\mathrel EH\wedge X\mathrel EY).     \tag{I10}
\]

\emph{Cartesian Product.}
Use Kuratowski pairs and the already constructed pairing graphs.  For
inputs with $n,m$ vertices, take disjoint copies of the inputs and, for
each $i<n,j<m$, add vertices $s_{ij},d_{ij},o_{ij}$ with edges
\[
 s_{ij}\longrightarrow i,\qquad
 d_{ij}\longrightarrow i,n+j,\qquad
 o_{ij}\longrightarrow s_{ij},d_{ij}.                   \tag{I11}
\]
Add a final root pointing to $o_{ij}$ precisely when $i$ and $j$ are root
successors in the respective input graphs.  The construction has
\[
    n+m+3nm+1
\]
vertices, a fixed polynomial in the input sizes.  In fact, if the input
padding exponents are $k\geq n^2$ and $\ell\geq m^2$, then the new vertex
count $N$ satisfies $N\leq5(k+\ell)$ and
$N^2\leq25(k+\ell)^2$; one further application of $\omega_1$ to a fixed
power of the two padding fields supplies the required matrix bound.  By
(I9), $s_{ij},d_{ij},o_{ij}$ represent respectively
$\{G[i]\}$, $\{G[i],H[j]\}$, and
$\langle G[i],H[j]\rangle$.  Hence the new root represents exactly
$G\times H$.

\emph{$\Delta_0$-Separation.}
Fix $\varphi(x,\bar v)\in\Delta_0$ and codes $H,\bar P$.  Copy $H$ and add
a root $u$ with
\[
 u\longrightarrow i
 \quad\longleftrightarrow\quad
 \operatorname{Edge}_H(r_H,i)\wedge
 \operatorname{Sat}_\varphi(H[i],\bar P).                \tag{I12}
\]
The predicate on the right is bounded by
Lemma~\ref{lem:local-bounded-truth}, so the graph-construction principle
produces a code $S_\varphi(H,\bar P)$.  For every $X\in D$,
\[
\begin{aligned}
 X\mathrel E S_\varphi(H,\bar P)
 &\longleftrightarrow X\mathrel EH
       \wedge\operatorname{Sat}_\varphi(X,\bar P)\\
 &\longleftrightarrow X\mathrel EH\wedge\varphi^I(X,\bar P).
\end{aligned}\tag{I13}
\]
Thus every translated Separation instance holds.  The domain is nonempty,
equality is a congruence, and all axioms in the chosen formulation of
$\mathsf{DB}_0$ have now been verified.
\end{proof}

\begin{corollary}\label{cor:mutual-interpretability}
In the usual definitionally equivalent languages,
\[
 \mathsf{DB}_0\equiv_{\mathrm{int}}\mathsf Q
 \equiv_{\mathrm{int}} I\Delta_0+\Omega_1,
 \qquad
 \mathsf{DB}_0\equiv_{\mathrm{int}}\mathsf{PA}^-.
\]
\end{corollary}
\begin{proof}
Proposition~\ref{prop:db0-interprets-q}, Wilkie's interpretation, and
Theorem~\ref{thm:dag-interpretation} give the cycle
\[
 \mathsf{DB}_0\triangleright\mathsf Q
 \triangleright I\Delta_0+\Omega_1
 \triangleright\mathsf{DB}_0.
\]
The usual discretely ordered semiring presentation of $\mathsf{PA}^-$
interprets $\mathsf Q$ by defining successor as $x+1$, while
$I\Delta_0+\Omega_1$ interprets $\mathsf{PA}^-$ by forgetting induction
and $\Omega_1$.  Composing these arrows in both directions gives the last
claim.
\end{proof}

\begin{remark}
Mutual interpretability is all that is asserted.  We have not shown that
the composites are definably isomorphic to the identity interpretations,
which would be a bi-interpretability result.  The use of DAGs is also
essential at this strength.  In Ackermann coding the singleton of a set
with code $a$ requires a bit at position $a$, hence a number at least
$2^a$.  The theory $S$ uniformly supplies only powers whose exponents are
fixed polynomials in input length.  DAG operations have precisely that
size; Power Set does not.  This distinguishes the construction from the
Ackermann interpretation over $I\Delta_0+\exp$ in \cite{Pe09}.
\end{remark}

\section{The Gaifman Splitting Property}

\begin{theorem}\label{thm:strong-splitting}
    Let $\mathcal M\vDash\mathsf{DB}_0$ and $n\geq0$.  The following are
    equivalent.
    \begin{enumerate}
        \item $\mathcal M$ has $\mathbf{end}_{n+1}$.
        \item $\mathcal M\vDash\mathsf{Coll}_s(\Sigma_{n+1})$.
        \item If $n=0$, then
        $\Mstar\preccurlyeq_{\Sigma_1}\mathcal N$ whenever
        $\mathcal M\preccurlyeq_{\Sigma_1}\mathcal N$ and
        $\mathcal N\vDash\mathsf{DB}_0$.  If $n\geq1$, then
        $\mathcal M\vDash\mathsf{Coll}(\Sigma_n)$ and
        $\Mstar\preccurlyeq_{\Sigma_{n+1}}\mathcal N$ whenever
        \[
           \mathcal M\preccurlyeq_{\Sigma_{n+1}}\mathcal N
           \quad\text{and}\quad
           \mathcal N\vDash
           \mathsf{DB}_0+\mathsf{Coll}(\Sigma_n).
        \]
    \end{enumerate}
\end{theorem}
\begin{proof}
    We argue by induction on $n$.

    $(1)\Rightarrow(2)$.  By
    Theorem~\ref{thm:strong-parameters} it is enough to prove the
    parameter-free scheme.  If it fails for a
    $\Sigma_{n+1}$ formula $\varphi(x,y)$, choose $p\in M$ such that
    \[
       \mathcal M\vDash
       \forall q\,\exists x\in p
       \bigl(\exists y\,\varphi(x,y)\wedge
             \forall y\in q\,\neg\varphi(x,y)\bigr).
       \tag{1}
    \]
    Add a constant $c$ and consider
    \[
    \begin{split}
       T={}&\ElDiag(\mathcal M)
       \cup\{c\in c_p,\ \exists y\,\varphi(c,y)\}\\
       &{}\cup\{\forall y\in c_q\,\neg\varphi(c,y):q\in M\}.
    \end{split}
    \]
    A finite subset mentions only $q_0,\ldots,q_{r-1}$.  Interpret $c$
    in $\mathcal M$ using (1) with
    $q=q_0\cup\cdots\cup q_{r-1}$.  Thus $T$ is finitely satisfiable.

    Let $\mathcal N\vDash T$, identifying the named copy of
    $\mathcal M$ with an elementary submodel.  Since $c\in p$, we have
    $c\in M^*$.  By $\mathbf{end}_{n+1}$, the assertion
    $\exists y\,\varphi(c,y)$ has a witness $b\in M^*$.  Choose
    $q\in M$ with $\mathcal N\vDash b\in q$.  Agreement between
    $\Mstar$ and $\mathcal N$ now contradicts the $q$-instance included
    in $T$.

    $(2)\Rightarrow(3)$.  The lower-level Collection assertion in (3)
    follows from Theorems~\ref{thm:collection-variants} and
    \ref{thm:strong-collection-equivalence}.  Let
    $\mathcal M\preccurlyeq_{\Sigma_{n+1}}\mathcal N$ satisfy the stated
    hypotheses.  The induction hypothesis gives
    $\Mstar\preccurlyeq_{\Sigma_n}\mathcal N$; for $n=0$ this is simply
    bounded absoluteness between an end submodel and the larger model.

    Suppose $a\in M^*$ and
    $\mathcal N\vDash\exists y\,\psi(a,y)$, where $\psi(x,y)$ is $\Pi_n$.
    Choose $p\in M$ with $\mathcal N\vDash a\in p$.  Strong Collection
    in $\mathcal M$ gives $q\in M$ such that
    \[
       \mathcal M\vDash\forall x\in p
       \bigl(\exists y\,\psi(x,y)\longrightarrow
              \exists y\in q\,\psi(x,y)\bigr).             \tag{2}
    \]
    We spell out the complexity of this assertion.  Its matrix is
    \[
       \forall y\,\neg\psi(x,y)
       \ \vee\ \exists y\in q\,\psi(x,y).                  \tag{2a}
    \]
    The first disjunct is $\Pi_{n+1}$.  If $n=0$, the second is bounded.
    If $n\geq1$, it is a $\Delta_0(\Pi_n)$ formula and is therefore
    $\Pi_n$ by Lemma~\ref{lem:bounded-contraction}; the required
    $\mathsf{Coll}(\Sigma_n)$ holds in both $\mathcal M$ and $\mathcal N$.
    Thus (2a), and then its bounded universal closure over $x\in p$, is
    $\Pi_{n+1}$ in both models.  Since
    $\mathcal M\preccurlyeq_{\Sigma_{n+1}}\mathcal N$, (2) holds in
    $\mathcal N$.

    Substituting $a$ in (2), the first disjunct is false by assumption, so
    $\mathcal N\vDash\exists y\in q\,\psi(a,y)$.  Every member of the old
    set $q$ belongs to $M^*$, and
    $\Mstar\preccurlyeq_{\Sigma_n}\mathcal N$ transfers the $\Pi_n$
    matrix to $\Mstar$.  Hence the required witness lies in $M^*$.  The
    $\Sigma_{n+1}$ Tarski--Vaught test proves
    $\Mstar\preccurlyeq_{\Sigma_{n+1}}\mathcal N$.

    $(3)\Rightarrow(1)$ follows by applying (3) to a full elementary
    extension.
\end{proof}

The next complexity observation will be used to transfer fragments of
Collection across partial elementary extensions.  The argument is due to
Zachiri McKenzie (personal communication).
\begin{theorem}\label{thm:collection-axiomatization}
    For $n\geq1$, each of
    $\mathsf{DB}_0+\mathsf{Coll}(\Sigma_n)$ and
    $\mathsf{DB}_0+\mathsf{Coll}_s(\Sigma_n)$ has an axiomatization by
    $\Pi_{n+2}$ sentences.
    At level $n=0$, both theories have $\Pi_3$ axiomatizations.
\end{theorem}
\begin{proof}
    The axioms of $\mathsf{DB}_0$ are $\Pi_2$.  We begin with a direct
    calculation for a $\Delta_0$ formula
    $\delta(x,y,\bar v)$.  Its ordinary Collection instance is equivalent
    to
    \[
    \forall\bar v\,\forall p\,\exists q\bigl[
       \exists x\in p\,\forall y\,\neg\delta(x,y,\bar v)
       \ \vee\
       \forall x\in p\,\exists y\in q\,\delta(x,y,\bar v)
       \bigr].                                             \tag{A1}
    \]
    The first disjunct in brackets is $\Pi_1$ and the second is bounded.
    Moving the unbounded quantifiers to the front therefore puts (A1) in
    $\Pi_3$.  A strong $\Delta_0$-Collection instance has the form
    \[
    \forall\bar v\,\forall p\,\exists q\,\forall x\in p
       \bigl(\forall y\,\neg\delta(x,y,\bar v)
       \ \vee\ \exists y\in q\,\delta(x,y,\bar v)\bigr),   \tag{A2}
    \]
    since the reverse implication is automatic.  The same prenex
    calculation makes (A2) $\Pi_3$.  This proves the separately stated
    $n=0$ bound.

    For $n=1$, Theorem~\ref{thm:collection-variants} replaces
    $\Sigma_1$-Collection by $\Pi_0$-Collection, so (A1) gives a
    $\Pi_3$ axiomatization.  By
    Theorem~\ref{thm:strong-collection-equivalence}, it remains, for strong
    Collection, to axiomatize $\Sigma_1$-Separation.  If
    $\sigma(x,\bar v)$ is $\Sigma_1$, its instance can be written
    \[
    \forall\bar v\,\forall p\,\exists q\bigl[
      \forall x\in p\,(\sigma(x,\bar v)\rightarrow x\in q)
      \ \wedge\
      \forall x\in q\,(x\in p\wedge\sigma(x,\bar v))
      \bigr].                                             \tag{A3}
    \]
    Relative to $\Sigma_1$-Collection, the first conjunct is $\Pi_1$ and,
    by Lemma~\ref{lem:bounded-contraction}, the second is $\Sigma_1$.
    Thus the bracket in (A3) can be prenexed with one existential and one
    universal block, and its universal closure is $\Pi_3$.

    Suppose inductively that the assertion has been proved through level
    $n\geq1$.  To treat level $n+1$, include the $\Pi_{n+2}$ axioms for the
    lower Collection theory and replace $\Sigma_{n+1}$-Collection by the
    equivalent $\Pi_n$-Collection scheme.  For
    $\varphi(x,y,\bar v)\in\Pi_n$, the corresponding version of (A1) has
    first disjunct
    \[
       \exists x\in p\,\forall y\,\neg\varphi(x,y,\bar v),
    \]
    which is $\Pi_{n+1}$ after bounded contraction, and second disjunct
    \[
       \forall x\in p\,\exists y\in q\,\varphi(x,y,\bar v),
    \]
    which is $\Pi_n$ by
    Lemma~\ref{lem:bounded-contraction}.  Pulling out $q$ and then taking
    the universal closure gives a $\Pi_{n+3}$ sentence.

    Finally let $\sigma\in\Sigma_{n+1}$ in (A3).  Its first inclusion is
    $\Pi_{n+1}$ and its second inclusion is $\Sigma_{n+1}$, again by
    bounded contraction over the lower theory.  Prenexing their conjunction
    after the existential quantifier for $q$ gives a $\Pi_{n+3}$ sentence.
    This proves the induction step for both ordinary and strong Collection.
\end{proof}

\begin{theorem}\label{thm:cofinal-transfer}
    Let $\mathcal M,\mathcal N\vDash\mathsf{DB}_0$ with
    $\mathcal M\preccurlyeq_{\Delta_0}\mathcal N$ and
    $\mathcal M\subseteq_{\mathsf{cf}}\mathcal N$.  For $n\geq0$:
    \begin{enumerate}
        \item if $\mathcal M\vDash\mathsf{Coll}(\Sigma_{n+1})$, then
        $\mathcal M\preccurlyeq_{\Sigma_{n+2}}\mathcal N$;
        \item if $\mathcal M\vDash\mathsf{Coll}_s(\Sigma_{n+1})$, then
        $\mathcal N\vDash\mathsf{Coll}_s(\Sigma_{n+1})$.
    \end{enumerate}
\end{theorem}
\begin{proof}
    We prove (1) and (2) simultaneously by induction on $n$.

    First consider (1) for $n=0$.  Cofinality and
    $\Delta_0$-elementarity already give
    $\mathcal M\preccurlyeq_{\Sigma_1}\mathcal N$: if
    $\mathcal N\vDash\exists x\,\delta(x,a)$ with $a\in M$ and
    $\delta\in\Delta_0$, choose a witness $b\in N$ and then $p\in M$ with
    $\mathcal N\vDash b\in p$.  The bounded assertion
    $\exists x\in p\,\delta(x,a)$ transfers to $\mathcal M$.

    Now suppose
    \[
       \mathcal N\vDash\exists y\,\forall z\,\theta(a,y,z),
       \qquad a\in M,\quad\theta\in\Delta_0,               \tag{T1}
    \]
    but (T1) fails in $\mathcal M$.  For every $p\in M$ we then have
    \[
       \mathcal M\vDash
       \forall y'\in p\,\exists z\,\neg\theta(a,y',z).
    \]
    By $\mathsf{Coll}(\Sigma_1)$, choose $q_p\in M$ such that
    \[
       \mathcal M\vDash
       \forall y'\in p\,\exists z\in q_p\,
          \neg\theta(a,y',z).                              \tag{T2}
    \]
    The last formula is bounded, so it also holds in $\mathcal N$.
    Add a constant $c$ and consider
    \[
       T=\ElDiag(\mathcal N)\cup
       \{\forall y'\in c_p\,\exists z\in c\,
          \neg\theta(c_a,y',z):p\in M\}.                  \tag{T3}
    \]
    A finite part of the second set mentions $p_0,\ldots,p_{k-1}$.
    Replace them by their union $p$ and interpret $c$ in $\mathcal N$ as
    the bound $q_p$ from (T2).  Hence (T3) is finitely satisfiable.

    Let $\mathcal K\vDash T$, and let $b\in N$ witness (T1).  Cofinality
    gives $p\in M$ with $\mathcal N\vDash b\in p$.  The $p$-sentence in
    (T3) yields in $\mathcal K$ a $z$ with
    $\neg\theta(a,b,z)$, whereas $\forall z\,\theta(a,b,z)$ belongs to
    $\ElDiag(\mathcal N)$.  This contradiction proves
    $\mathcal M\preccurlyeq_{\Sigma_2}\mathcal N$.

    We next prove (2) at level $n$, assuming (1) at that level.  By
    Theorem~\ref{thm:strong-parameters}, it is enough to transfer a
    parameter-free strong-Collection instance.  Let
    $\varphi(x,y)\in\Sigma_{n+1}$ and $p\in N$.  Choose $P\in M$ with
    $\mathcal N\vDash p\in P$.  Strong Collection in $\mathcal M$ supplies
    $q\in M$ such that
    \[
      \mathcal M\vDash\forall x\in\bigcup P
       \bigl(\exists y\,\varphi(x,y)\longrightarrow
             \exists y\in q\,\varphi(x,y)\bigr).           \tag{T4}
    \]
    To see the complexity without assuming Collection in $\mathcal N$,
    negate (T4).  Its negation is
    \[
       \exists x\in\bigcup P\bigl(
          \exists y\,\varphi(x,y)\ \wedge\
          \forall y\in q\,\neg\varphi(x,y)
       \bigr).                                             \tag{T4a}
    \]
    The two conjuncts are respectively $\Sigma_{n+1}$ and
    $\Delta_0(\Pi_{n+1})$; their conjunction has a $\Sigma_{n+2}$ prenex
    form.  Hence (T4) is $\Pi_{n+2}$.  This calculation uses no Collection
    in $\mathcal N$.  Part (1) gives
    $\mathcal M\preccurlyeq_{\Sigma_{n+2}}\mathcal N$, so (T4) transfers
    to $\mathcal N$.  Since $p\in P$, every member of $p$ belongs to
    $\bigcup P$, and the same $q$ proves the desired strong-Collection
    instance in $\mathcal N$.

    Finally let $n\geq1$ and assume both assertions below level $n$.  From
    $\mathsf{Coll}(\Sigma_{n+1})$ in $\mathcal M$ and
    Theorem~\ref{thm:weak-to-strong}, we obtain
    $\mathsf{Coll}_s(\Sigma_n)$ in $\mathcal M$; the induction hypothesis
    for (2) transfers it to $\mathcal N$.  The induction hypothesis for
    (1) also gives
    $\mathcal M\preccurlyeq_{\Sigma_{n+1}}\mathcal N$.

    Suppose
    \[
       \mathcal N\vDash\exists x\,\forall y\,\theta(x,y,a),
       \qquad a\in M,\quad\theta\in\Sigma_n.               \tag{T5}
    \]
    Cofinality puts a witness in some $p\in M$.  Hence, for every $q\in M$,
    \[
       \mathcal N\vDash
       \exists x\in p\,\forall y\in q\,\theta(x,y,a).     \tag{T6}
    \]
    Both models satisfy $\mathsf{Coll}(\Sigma_n)$.  By
    Lemma~\ref{lem:bounded-contraction}, (T6) has the same $\Sigma_n$
    equivalent in both models, so lower-level elementarity puts (T6) in
    $\mathcal M$ for every $q\in M$.

    If no $x\in p$ satisfied $\forall y\,\theta(x,y,a)$ in $\mathcal M$,
    then
    \[
       \mathcal M\vDash
       \forall x\in p\,\exists y\,\neg\theta(x,y,a).
    \]
    Here $\neg\theta\in\Pi_n$, and
    $\mathsf{Coll}(\Pi_n)$ is equivalent to
    $\mathsf{Coll}(\Sigma_{n+1})$.  It would therefore produce one
    $q\in M$ containing a counterexample for every $x\in p$, contradicting
    the $q$-instance of (T6).  Thus (T5) reflects to $\mathcal M$, completing
    the induction for (1).
\end{proof}
\begin{remark}
    The last part can alternatively be proved by the compactness argument
    used in the base case; compare \cite{Ku24}.
\end{remark}

We need one closure observation before sharpening
Theorem~\ref{thm:strong-splitting}.
\begin{lemma}\label{lem:membership-hull-db0}
    If $\mathcal M\preccurlyeq\mathcal N$ and
    $\mathcal M\vDash\mathsf{Coll}_s(\Sigma_1)$, then its membership hull
    $\Mstar$ is a model of $\mathsf{DB}_0$.
\end{lemma}
\begin{proof}
    Each operation occurring in the axioms of $\mathsf{DB}_0$, including
    each fixed $\Delta_0$ separator, is uniquely specified by a
    $\Sigma_1$ formula.  Put the finitely many arguments in old sets and
    apply strong Collection in $\mathcal M$ to all possible tuples from
    those sets.  The resulting old bound contains the value of the
    operation on every tuple from $M^*$.  Hence $M^*$ is closed under
    Empty Set, Pairing, Union, Cartesian Product, and every required
    $\Delta_0$ separator.  Extensionality is inherited from $\mathcal N$.
\end{proof}

The following refinement lowers the auxiliary Collection fragment required
of the larger model.
\begin{theorem}\label{thm:refined-strong-splitting}
    Let $\mathcal M\vDash\mathsf{DB}_0$ and $n\geq0$.  The following are
    equivalent.
    \begin{enumerate}
        \item $\mathcal M\vDash\mathsf{Coll}_s(\Sigma_{n+1})$.
        \item For $n=0$, the extension property in
        Theorem~\ref{thm:strong-splitting}(3) holds.
        For $n=1$, $\mathcal M\vDash\mathsf{Coll}(\Sigma_1)$ and
        $\Mstar\preccurlyeq_{\Sigma_2}\mathcal N$ whenever
        $\mathcal M\preccurlyeq_{\Sigma_2}\mathcal N$ and
        $\mathcal N\vDash\mathsf{DB}_0$.
        For $n\geq2$, $\mathcal M\vDash\mathsf{Coll}(\Sigma_n)$ and
        $\Mstar\preccurlyeq_{\Sigma_{n+1}}\mathcal N$ whenever
        $\mathcal M\preccurlyeq_{\Sigma_{n+1}}\mathcal N$ and
        $\mathcal N\vDash
          \mathsf{DB}_0+\mathsf{Coll}(\Sigma_{n-1})$.
        \item For every
        $\mathcal M\preccurlyeq_{\Sigma_{n+1}}\mathcal N$ with
        $\mathcal N\vDash\mathsf{DB}_0$, one has
        $\Mstar\preccurlyeq_{\Sigma_{n+1}}\mathcal N$.
    \end{enumerate}
\end{theorem}
\begin{proof}
    Clause (3) implies (1) by Theorem~\ref{thm:strong-splitting}.
    Clause (3) trivially implies (2), so it remains to prove
    $(1)\Rightarrow(3)$.  We do this by induction on $n$.  The case $n=0$
    is Theorem~\ref{thm:strong-splitting}.

    Suppose $n\geq1$ and
    $\mathcal M\preccurlyeq_{\Sigma_{n+1}}\mathcal N$ with
    $\mathcal N\vDash\mathsf{DB}_0$.  Strong Collection in $\mathcal M$
    implies $\mathsf{Coll}(\Sigma_n)$.  If $n\geq2$, then
    $\mathsf{DB}_0+\mathsf{Coll}(\Sigma_{n-1})$ has a
    $\Pi_{n+1}$ axiomatization by
    Theorem~\ref{thm:collection-axiomatization}; since it holds in
    $\mathcal M$, $\Sigma_{n+1}$-elementarity transfers it to
    $\mathcal N$.  When $n=1$ no such auxiliary Collection is needed in
    the calculation below.

    By the induction hypothesis,
    $\Mstar\preccurlyeq_{\Sigma_n}\mathcal N$.  To perform the next
    Tarski--Vaught step, write the relevant $\Pi_n$ matrix as
    \[
       \psi(x,y):=\forall z\,\theta(x,y,z),
       \qquad \theta\in\Sigma_{n-1}.
    \]
    Suppose $a\in M^*$ and
    $\mathcal N\vDash\exists y\,\psi(a,y)$.  Choose $p\in M$ with
    $a\in^{\mathcal N}p$.  Strong Collection supplies $q\in M$ such that
    \[
      \mathcal M\vDash\forall x\in p\left(
        \forall y\,\exists z\,\neg\theta(x,y,z)
        \ \vee\ \exists y\in q\,\forall z\,\theta(x,y,z)
      \right).                                             \tag{3a}
    \]
    Instead of transferring (3a), transfer its following consequence:
    \[
    \forall x\in p\left(
       \forall y\,\exists z\,\neg\theta(x,y,z)
       \ \vee\
       \forall w\,\exists y\in q\,\forall z\in w\,
          \theta(x,y,z)
       \right).                         \tag{3}
    \]

    Here is the promised complexity calculation.  The first disjunct is
    $\Pi_{n+1}$.  In the second disjunct, if $n=1$ then
    \[
       \exists y\in q\,\forall z\in w\,\theta(x,y,z)
    \]
    is bounded, so the unbounded quantifier $\forall w$ makes the second
    disjunct $\Pi_1$.  If $n\geq2$, its matrix is a
    $\Delta_0(\Sigma_{n-1})$ formula.  The auxiliary
    $\mathsf{Coll}(\Sigma_{n-1})$ and
    Lemma~\ref{lem:bounded-contraction} replace it uniformly by a
    $\Sigma_{n-1}$ formula; the leading $\forall w$ therefore makes the
    second disjunct $\Pi_n$.  In either case the disjunction, followed by
    the bounded universal quantifier $x\in p$, is $\Pi_{n+1}$.  Thus (3)
    transfers from $\mathcal M$ to $\mathcal N$.

    Substitute the chosen $a$.  Its $\mathcal N$-witness to
    $\exists y\,\forall z\,\theta(a,y,z)$ makes the first disjunct false,
    so
    \[
       \mathcal N\vDash
       \forall w\,\exists y\in q\,\forall z\in w\,
          \theta(a,y,z).                                  \tag{3b}
    \]
    Formula (3b) is $\Pi_n$ by the calculation just given.  Hence
    $\Mstar\preccurlyeq_{\Sigma_n}\mathcal N$ puts (3b) in $\Mstar$.

    By Lemma~\ref{lem:membership-hull-db0},
    $\Mstar\vDash\mathsf{DB}_0$.  Theorem~\ref{thm:cofinal-transfer},
    applied to
    $\mathcal M\subseteq_{\mathsf{cf}}\Mstar$, transfers strong
    Collection to $\Mstar$.  In particular it has
    $\mathsf{Coll}(\Sigma_n)$, equivalently
    $\mathsf{Coll}(\Pi_{n-1})$.  If no $y\in q$ satisfied
    $\forall z\,\theta(a,y,z)$ in $\Mstar$, this Collection instance would
    put, into one set $w$, a counterexample $z$ for every $y\in q$.  That
    contradicts the $w$-instance of (3b).  Therefore
    \[
       \Mstar\vDash\exists y\in q\,\forall z\,\theta(a,y,z),
    \]
    supplying the required witness in $M^*$.  This completes the
    $\Sigma_{n+1}$ Tarski--Vaught step and the induction.
\end{proof}

\begin{theorem}[Refinement of Kurahashi--Minami]
\label{thm:arithmetic-splitting}
    Let $\mathcal M\vDash\mathsf{PA}^-$ and $n\geq0$.  Then
    $\mathcal M\vDash\mathsf{Coll}_s(\Sigma_{n+1})$ if and only if
    \[
       \operatorname{sup}_{\mathcal N}(M)
       \preccurlyeq_{\Sigma_{n+1}}\mathcal N
    \]
    for every
    $\mathcal M\preccurlyeq_{\Sigma_{n+1}}\mathcal N$ with
    $\mathcal N\vDash\mathsf{PA}^-$.
\end{theorem}
\begin{proof}
    Put $H=\operatorname{sup}_{\mathcal N}(M)$.  The arithmetic operations
    make $H$ a $\mathsf{PA}^-$-substructure of $\mathcal N$.  We prove the
    forward implication simultaneously for all $n$, so at level $n+1$ the
    induction hypothesis already gives
    $H\preccurlyeq_{\Sigma_n}\mathcal N$.

    Write a $\Pi_n$ matrix as
    \[
       \psi(x,y)=\forall z\,\theta(x,y,z),
       \qquad\theta\in\Sigma_{n-1},\quad n\geq1.
    \]
    Let $a\in H$, choose $p\in M$ with $a<p$ in $\mathcal N$, and suppose
    $\mathcal N\vDash\exists y\,\psi(a,y)$.  Strong Collection in
    $\mathcal M$ gives a bound $q\in M$ and, exactly as in (3a), implies
    \[
    \forall x<p\left(
       \forall y\,\exists z\,\neg\theta(x,y,z)
       \ \vee\
       \forall w\,\exists y<q\,\forall z<w\,
          \theta(x,y,z)
       \right).                                           \tag{3c}
    \]
    The first disjunct is $\Pi_{n+1}$.  In the second, the quantifiers
    $y<q$ and $z<w$ are bounded arithmetic quantifiers, so its matrix is
    still $\Sigma_{n-1}$ and the leading $\forall w$ makes it $\Pi_n$.
    The outer quantifier $x<p$ is also bounded.  Consequently (3c) is a
    $\Pi_{n+1}$ formula, with no auxiliary Collection assumption on
    $\mathcal N$.  It therefore transfers to $\mathcal N$.

    At $x=a$ the first disjunct is false, and the second disjunct, being
    $\Pi_n$, transfers to $H$.  The arithmetic cofinal-transfer argument
    used in \cite[Theorem~3.1]{Ku24} gives $H$ the required lower
    $\mathsf{Coll}(\Sigma_n)$ instance.  If no $y<q$ worked for all $z$,
    that instance would bound one counterexample for every $y<q$ below a
    single $w$, contradicting (3c).  Hence $H$ contains a witness, and the
    $\Sigma_{n+1}$ Tarski--Vaught test is complete.  When $n=0$, the same
    proof omits the displayed unbounded $z$-block and is just bounded
    absoluteness.

    Conversely, suppose the asserted splitting property holds for every
    $\Sigma_{n+1}$-elementary extension.  If parameter-free strong
    $\Sigma_{n+1}$-Collection failed for $\varphi(x,y)$ at a bound $p$,
    compactness applied to
    \[
       \ElDiag(\mathcal M)\cup\{c<p,\ \exists y\,\varphi(c,y)\}
       \cup\{\forall y<q\,\neg\varphi(c,y):q\in M\}
    \]
    would give an elementary extension $\mathcal N$ with $c\in H$.
    Elementarity of $H$ would provide a witness $b\in H$; some $q\in M$
    bounds $b$, contradicting the corresponding sentence in the type.
    Parameter coding then gives full strong Collection.  This proves the
    stated refinement of \cite[Theorem~3.1]{Ku24}.
\end{proof}

\begin{theorem}\label{thm:ordinary-collection}
    Let $\mathcal M\vDash\mathsf{DB}_0$ and $n\geq0$.  The following are
    equivalent.
    \begin{enumerate}
        \item $\mathcal M\vDash\mathsf{Coll}(\Sigma_{n+1})$.
        \item $\mathcal M$ has $\mathbf{end}_n$ and
        $\mathbf{COF}_{n+2}$.
        \item $\mathcal M$ has $\mathbf{end}_n$ and
        $\mathbf{cof}_{n+2}$.
    \end{enumerate}
\end{theorem}
\begin{proof}
    $(1)\Rightarrow(2)$ follows from
    Theorems~\ref{thm:strong-splitting} and
    \ref{thm:cofinal-transfer} (with the
    $n=0$ end assertion given by bounded absoluteness), and
    $(2)\Rightarrow(3)$ is immediate.  We prove $(3)\Rightarrow(1)$.

    By Theorem~\ref{thm:collection-variants} it is enough to collect
    witnesses for a $\Pi_n$ formula
    $\varphi(x,y,v)$.  Suppose $p,v\in M$ and
    \[
       \mathcal M\vDash
       \forall q\,\exists x\in p\,\forall y\in q\,
       \neg\varphi(x,y,v).                               \tag{4}
    \]
    Add a constant $c$ and let $T$ consist of $\ElDiag(\mathcal M)$ and,
    for every $a\in M$, the sentences
    \[
       c_a\in c,\qquad
       \forall z\in c_a\ (z\in c),qquad
       c\notin c_a.                                      \tag{5}
    \]
    We verify finite satisfiability, since the last part of (5) was missing
    from the informal compactness argument.  Given finitely many
    requirements, first form an old set $b$ containing the named elements
    and every member of the finitely many named sets.  Let $q$ be the union
    of the finitely many sets from which $c$ must be excluded.  Some
    $d\in M$ has $b\cup\{d\}\notin q$: otherwise $\bigcup q$ would contain
    every $d\in M$, contrary to the Russell argument used in
    Theorem~\ref{thm:db0-exact-cover}.
    Interpreting $c$ as $b\cup\{d\}$ satisfies the finite fragment.

    Let $\mathcal N\vDash T$.  The last group in (5) gives
    $c\notin M^*$, while the first two groups give $M^*\subseteq
    \operatorname{Ext}_{\mathcal N}(c)$.  Transfer (4) to $\mathcal N$
    and substitute $c$ for $q$.  We obtain $x\in p$, hence $x\in M^*$,
    such that
    $\mathcal N\vDash\neg\varphi(x,y,v)$ for every $y\in M^*$.
    Since $\neg\varphi$ is $\Sigma_n$, $\mathbf{end}_n$ yields
    \[
       \Mstar\vDash
       \exists x\in p\,\forall y\,\neg\varphi(x,y,v).
    \]
    This is a $\Pi_{n+1}$, hence $\Sigma_{n+2}$, assertion.
    Property $\mathbf{cof}_{n+2}$ reflects it to $\mathcal M$.  Thus (4)
    always yields an $x\in p$ with no witness at all, which is the
    contrapositive form of Collection.
\end{proof}

\begin{theorem}\label{thm:end-implies-cof}
    For every $n\geq0$, property $\mathbf{end}_n$ implies both
    $\mathbf{cof}_{n+1}$ and $\mathbf{COF}_{n+1}$.
\end{theorem}
\begin{proof}
    We first prove $\mathbf{cof}_{n+1}$ by induction on $n$.  For $n=0$,
    let $\mathcal M\preccurlyeq\mathcal N$ and suppose a bounded formula
    $\delta(a,b)$, with $a\in M$ and $b\in M^*$, witnesses a
    $\Sigma_1$ assertion in $\Mstar$.  Property $\mathbf{end}_0$ is bounded
    absoluteness between $\Mstar$ and $\mathcal N$, so the same assertion
    holds in $\mathcal N$ and hence in $\mathcal M$.  The reverse direction
    uses the witness from $M$.  Thus
    $\mathcal M\preccurlyeq_{\Sigma_1}\Mstar$.

    Now suppose $n\geq1$ and let
    $\mathcal M\preccurlyeq\mathcal N$.  Since
    $\mathbf{end}_n$ implies $\mathbf{end}_{n-1}$, the induction hypothesis
    and the definition of $\mathbf{end}_n$ give
    \[
       \mathcal M\preccurlyeq_{\Sigma_n}\Mstar
       \preccurlyeq_{\Sigma_n}\mathcal N.
    \]
    Let $\exists x\,\psi(x,a)$ be $\Sigma_{n+1}$ with
    $\psi\in\Pi_n$ and $a\in M$.  If it holds in $\Mstar$, its witness is
    in $M^*$ and $\Sigma_n/\Pi_n$ agreement puts the same instance in
    $\mathcal N$; full elementarity then reflects the existential assertion
    to $\mathcal M$.  Conversely, a witness from $M$ satisfies the
    $\Pi_n$ matrix in $\Mstar$ by the left-hand agreement above.  Hence
    $\mathcal M\preccurlyeq_{\Sigma_{n+1}}\Mstar$, proving
    $\mathbf{cof}_{n+1}$.

    It remains to obtain the uniform property $\mathbf{COF}_{n+1}$.  For
    $n=0$, every $\Delta_0$-elementary cofinal extension is
    $\Sigma_1$-elementary by the bounded-witness argument used at the start
    of Theorem~\ref{thm:cofinal-transfer}.  For $n\geq1$, we have just proved
    $\mathbf{cof}_{n+1}$, while $\mathbf{end}_n$ implies
    $\mathbf{end}_{n-1}$.  Theorem~\ref{thm:ordinary-collection}, applied
    one level lower, gives
    $\mathcal M\vDash\mathsf{Coll}(\Sigma_n)$.  Part (1) of
    Theorem~\ref{thm:cofinal-transfer} now makes every
    $\Delta_0$-elementary cofinal $\mathsf{DB}_0$-extension of
    $\mathcal M$ a $\Sigma_{n+1}$-elementary extension.  This is
    $\mathbf{COF}_{n+1}$.
\end{proof}

\begin{theorem}\label{thm:full-splitting}
    For $\mathcal M\vDash\mathsf{DB}_0$, the following are equivalent.
    \begin{enumerate}
        \item $\mathcal M\vDash\mathsf{Coll}$.
        \item For every $\mathcal M\preccurlyeq\mathcal N$,
        $\mathcal M\preccurlyeq\Mstar\preccurlyeq\mathcal N$.
        \item For every $\mathcal M\preccurlyeq\mathcal N$,
        $\Mstar\preccurlyeq\mathcal N$.
    \end{enumerate}
\end{theorem}
\begin{proof}
    Apply Theorems~\ref{thm:ordinary-collection} and
    \ref{thm:end-implies-cof} at every finite level.  Conversely, clause
    (3) gives every property $\mathbf{end}_n$, and
    Theorem~\ref{thm:ordinary-collection} then gives every finite fragment
    of Collection.
\end{proof}

\section{The Keisler--Morley Extension Property}

Kaufmann \cite{Ka81} refined the Keisler--Morley theorem by relating
partial elementary end extensions to fragments of Collection.  Keisler and
Morley's construction \cite{Ke68}, by contrast, uses rank stages and
partition relations.  We first isolate the conclusions requiring no
resolution and then turn to the stronger resolved setting.

\subsection{End extensions without resolutions}

\begin{theorem}\label{thm:taller-to-taller-star}
    If $\mathcal M\vDash\mathsf{DB}_0+\mathsf{TCo}$, every taller
    $\Sigma_2$-elementary extension of $\mathcal M$ is taller*.
\end{theorem}
\begin{proof}
    Let $c$ witness that $\mathcal N$ is taller than $\mathcal M$.
    Since $\mathsf{TCo}$ is $\Pi_2$ and
    $\mathcal M\preccurlyeq_{\Sigma_2}\mathcal N$, the model $\mathcal N$
    satisfies $\mathsf{TCo}$.  Choose an $\mathcal N$-transitive $t$
    with $c\in t$.  If $a\in M$, then
    $\mathcal{N}\vDash a\in c\in t$, and transitivity gives
    $\mathcal{N}\vDash a\in t$.  Thus $t$ witnesses taller*.
\end{proof}
\begin{remark}\label{rem:tco-and-tcl}
    The hypothesis $\mathsf{TCo}$ is substantive and will always be
    displayed explicitly.  In particular, Mathias constructs a model of
    Zermelo set theory with full Foundation in which $\mathsf{TCo}$ fails
    \cite[Section~12]{Ma06}.

    There is also a precise finite-chain implication in the other direction.
    Mathias calls a class term semi-suitable when quantification bounded by
    that term preserves bounded complexity, and records that $\omega$ is
    semi-suitable over his base theory with Set Foundation
    \cite[Definition~1.4 and Remark~1.6]{Ma06}.  In that setting put
    \[
    y\mathrel{\triangleleft}x\quad\Longleftrightarrow\quad
    \exists f\,\exists n\in\omega\biggl(
       \begin{array}{l}
       \operatorname{Fun}(f)\ \wedge\ \operatorname{dom}(f)=n+2\ \wedge\
       f(0)=x\ \wedge\\[-1mm]
       f(n+1)=y\ \wedge\
       \forall i\in n+1\;f(i+1)\in f(i)
       \end{array}
    \biggr).                                                \tag{TC}
    \]
    Here the usual finite-sequence coding is available in $\mathsf{DB}_0$.
    Semi-suitability of $\omega$ shows that the matrix following the single
    unbounded existential quantifier over $f$ can be taken bounded; hence
    $y\triangleleft x$ is $\Sigma_1$.

    Given $x$, use $\mathsf{TCo}$ to choose a transitive $u$ with $x\in u$,
    and use $\mathsf{Sep}(\Sigma_1)$ to form
    \[
       t=\{y\in u:y\triangleleft x\}.
    \]
    Every member of $x$ belongs to $t$, by a chain of length one.  Appending
    one value to a chain proves that $t$ is transitive.  Finally, induction
    along the finite domain of a chain shows that every transitive set
    containing $x$ as a subset contains $t$.  Thus $t$ is the transitive
    closure of $x$.  Consequently
    \[
      \mathsf{DB}_0+\mathsf{Set\ Foundation}+\mathsf{TCo}
        +\mathsf{Sep}(\Sigma_1)\ \vdash\ \mathsf{TCl}.
    \]
    This is the Separation version of the construction in
    \cite[Proposition~8.17]{Ma06}.  The Set Foundation qualification is
    essential to the cited semi-suitability calculation, so we do not claim
    the displayed implication over bare $\mathsf{DB}_0$.
\end{remark}

The first implication from extensions to Collection uses no resolution.
\begin{theorem}\label{thm:taller-implies-collection}
    Let $n\geq0$.
    If $\mathcal M\vDash
    \mathsf{DB}_0+\mathsf{Coll}(\Sigma_1)$ has a
    $\Sigma_{n+2}$-elementary taller end extension, then
    $\mathcal M\vDash\mathsf{Coll}(\Sigma_{n+2})$.
\end{theorem}
\begin{proof}
    We argue by induction on $n$.  The standing assumption gives the lower
    Collection fragment when $n=0$; at a later stage the induction
    hypothesis gives $\mathsf{Coll}(\Sigma_{n+1})$.

    Let $\mathcal M\preccurlyeq_{\Sigma_{n+2}}\mathcal N$ be an end
    extension, and choose $c\in N$ with
    $\mathcal{N}\vDash a\in c$ for every $a\in M$.  By
    Theorem~\ref{thm:collection-variants}, it is enough to prove Collection
    for a $\Pi_{n+1}$ formula, which we write as
    \[
       \forall z\,\varphi(x,y,z,\bar v),
       \qquad\varphi\in\Sigma_n.
    \]
    Suppose
    \[
       \mathcal{M}\vDash\forall x\in p\,\exists y\,
       \forall z\,\varphi(x,y,z,\bar v),
       \qquad p,\bar v\in M.
                                                               \tag{L1}
    \]
    End extension ensures that every $\mathcal{N}$-member of $p$ is in
    $M$.  For each such $x$, choose in $M$ a witness $y$ to the displayed
    formula.  Agreement for $\Pi_{n+1}$ formulas puts the same witness
    statement in $\mathcal N$, and taller-ness gives $y\in^{\mathcal N}c$.
    Therefore
    \[
       \mathcal N\vDash\exists q\,\forall x\in p\,\forall w\,
       \exists y\in q\,\forall z\in w\,
          \varphi(x,y,z,\bar v),                            \tag{L2}
    \]
    with $q=c$.

    We verify rather than assume the complexity of (L2).  If $n=0$, all
    quantifiers after $\forall w$ are bounded and its matrix is
    $\Delta_0$; hence (L2) is $\Sigma_2$.  Suppose $n\geq1$ and write
    \[
       \varphi(x,y,z,\bar v)\equiv
       \exists t\,\chi(x,y,z,t,\bar v),
       \qquad\chi\in\Pi_{n-1}.
    \]
    By induction $\mathcal M\vDash\mathsf{Coll}(\Sigma_n)$.  This theory
    has a $\Pi_{n+2}$ axiomatization by
    Theorem~\ref{thm:collection-axiomatization}, so it transfers to
    $\mathcal N$.  In either model Collection for the $\Pi_{n-1}$ matrix
    gives
    \begin{equation}
    \tag{L3}
    \begin{aligned}
      &\forall x\in p\,\forall w\,\exists y\in q\,
        \forall z\in w\,\exists t\,\chi(x,y,z,t,\bar v)\\
      &\quad\Longleftrightarrow\quad
      \forall x\in p\,\forall w\,\exists T\,\exists y\in q\,
        \forall z\in w\,\exists t\in T\,
           \chi(x,y,z,t,\bar v).
    \end{aligned}
    \end{equation}
    For right-to-left no Collection is needed; left-to-right collects the
    $t$-witnesses after $y$ has been chosen.  The matrix following
    $\exists T$ is a $\Delta_0(\Pi_{n-1})$ formula and contracts to
    $\Pi_{n-1}$.  Thus the right side of (L3) is $\Pi_{n+1}$, and the
    leading existential quantifier for $q$ makes (L2) $\Sigma_{n+2}$.
    Hence (L2) reflects to $\mathcal M$.

    Fix $x\in^{\mathcal M}p$.  If no $y\in q$ satisfied
    $\forall z\,\varphi(x,y,z,\bar v)$, then
    \[
       \forall y\in q\,\exists z\,
          \neg\varphi(x,y,z,\bar v).                       \tag{L4}
    \]
    The matrix in (L4) is $\Pi_n$.  Since
    $\mathsf{Coll}(\Pi_n)$ is equivalent to
    $\mathsf{Coll}(\Sigma_{n+1})$, the lower Collection fragment supplied
    by the induction hypothesis bounds all these counterexamples in one
    set $w$.  This contradicts the $x,w$ instance of (L2).  Consequently
    \[
       \mathcal M\vDash\exists q\,\forall x\in p\,
       \exists y\in q\,\forall z\,\varphi(x,y,z,\bar v),
    \]
    proving $\mathsf{Coll}(\Pi_{n+1})$, equivalently
    $\mathsf{Coll}(\Sigma_{n+2})$.
\end{proof}

\begin{theorem}\label{thm:kaufmann-taller-extension}
    Let $n\geq0$.  Every countable model
    $\mathcal M\vDash\mathsf{DB}_0+\mathsf{TCo}+
    \mathsf{Coll}(\Sigma_{n+2})$ has a
    $\Sigma_{n+2}$-elementary taller* end extension which satisfies
    $\mathsf{DB}_0+\mathsf{TCo}$.  If
    $\mathcal M\vDash\mathsf{Coll}$, the extension may be chosen
    elementary.
\end{theorem}
\begin{proof}
    This is Kaufmann's construction for the implication
    (ii)$\Rightarrow$(i) in his Theorem~1 \cite{Ka81}, with one extra
    constant.
    We describe the modification and refer to his Lemma~4 for the
    stratified Henkin construction.

    Expand the $\Pi_{n+2}$ diagram of $\mathcal M$ by a constant $c$,
    the sentences $c_a\in c$ for $a\in M$, and $\Trans(c)$.  Every finite
    part is realized in $\mathcal M$: Pairing collects the finitely many
    named elements and $\mathsf{TCo}$ puts that set inside a transitive
    set.  For each $a\in M$, use Kaufmann's end-extension type
    \[
       \Gamma_a(x)=\{x\in c_a\}\cup
       \{x\ne c_b:b\in\operatorname{Ext}_{\mathcal M}(a)\}.       \tag{6}
    \]
    Omitting all types (6) is exactly the assertion that no old set gains
    a new member.

    The local omission condition in \cite[Lemma~4]{Ka81} is verified in
    his proof of Theorem~1 by applying Collection to the family of local
    obstructions indexed by $b\in^{\mathcal M}a$.  Adding $c$ changes
    only the preliminary consistency test: a finite list of requirements
    on $c$ is absorbed into one set $d$, after which
    $\mathsf{TCo}$ supplies a transitive $t$ with $d\in t$.  The same
    $\mathsf{Coll}(\Sigma_{n+2})$ instance collects the obstructions, so
    Kaufmann's argument applies without further change.

    The resulting model $\mathcal N$ is a
    $\Sigma_{n+2}$-elementary end extension of $\mathcal M$.  The
    $\Pi_2$ axioms of $\mathsf{DB}_0+\mathsf{TCo}$ belong to the diagram,
    so $\mathcal N$ satisfies them.  Its
    interpretation of $c$ is transitive and contains every old element,
    hence is a taller* witness.  It is proper, since the Russell argument
    shows that no element of $M$ contains all of $M$.  With full
    Collection, Kaufmann's construction is run through all finite
    complexities, yielding an elementary end extension.
\end{proof}

\begin{theorem}\label{thm:countable-end-characterization}
    Let $n\geq0$, and let $\mathcal M$ be a countable model of
    $\mathsf{DB}_0+\mathsf{TCo}+\mathsf{Coll}(\Sigma_1)$.  The following
    are equivalent:
    \begin{enumerate}
        \item $\mathcal M$ has a $\Sigma_{n+2}$-elementary taller* end
        extension;
        \item $\mathcal M$ has a $\Sigma_{n+2}$-elementary taller end
        extension;
        \item $\mathcal M\vDash\mathsf{Coll}(\Sigma_{n+2})$.
    \end{enumerate}
    Moreover, $\mathcal M$ has an elementary taller* end extension if and
    only if $\mathcal M\vDash\mathsf{Coll}$.
\end{theorem}
\begin{proof}
    The implications (1)$\Rightarrow$(2)$\Rightarrow$(3) follow from
    Theorems~\ref{thm:taller-to-taller-star} and
    \ref{thm:taller-implies-collection}.
    Theorem~\ref{thm:kaufmann-taller-extension} gives
    (3)$\Rightarrow$(1).  Applying the same argument at every finite level
    proves the last assertion.
\end{proof}

The countability hypothesis cannot simply be removed.  Enayat's construction
gives the following sharp obstruction.

\begin{theorem}[{\cite[Theorem~5.3 and Proposition~5.4]{En24}}]
\label{thm:aleph-one-dead-end}
Every countable model of $\mathsf{ZF}$ has an elementary end extension to
an $\aleph_1$-like model which has no taller extension satisfying
$\mathsf{ZF}$.
\end{theorem}

For an $\aleph_1$-like model, however, the part of the universe relevant to
one Collection instance is countable.  This yields a local version of the
preceding characterization without assuming any resolution.

\begin{theorem}\label{thm:aleph-one-local-characterization}
Let $n\geq0$, and let $\mathcal M$ be an $\aleph_1$-like model of
\[
   \mathsf{DB}_0+\mathsf{TCo}+\mathsf{Coll}(\Sigma_1).
\]
The following are equivalent.
\begin{enumerate}
    \item $\mathcal M\vDash\mathsf{Coll}(\Sigma_{n+2})$.
    \item For every countable $U\subseteq M$ there is
    $\mathcal N\vDash\mathsf{DB}_0+\mathsf{TCo}$ such that:
    \begin{enumerate}
        \item $U\subsetneq N$ and
        $(\mathcal N,u)_{u\in U}\equiv_{\Sigma_{n+2}}
        (\mathcal M,u)_{u\in U}$;
        \item every $u\in U$ is fixed:
        $\operatorname{Ext}_{\mathcal N}(u)=
        \operatorname{Ext}_{\mathcal M}(u)$;
        \item some $c\in N$ satisfies
        $\mathcal N\vDash u\in c$ for every $u\in U$.
    \end{enumerate}
\end{enumerate}
With full Collection, the corresponding assertion uses full elementary
equivalence.
\end{theorem}
\begin{proof}
Assume (1).  Starting from $U$, form its external downward closure
\[
 U_0=U,\qquad
 U_{k+1}=U_k\cup
 \bigcup_{u\in U_k}\operatorname{Ext}_{\mathcal M}(u),
 \qquad
 \widehat U=\bigcup_{k<\omega}U_k.                         \tag{E1}
\]
    Every stage is countable, since each
    $\operatorname{Ext}_{\mathcal M}(u)$ is countable and a countable union
    of countable sets is countable.  We now construct the required
    elementary submodel directly.  In the external metatheory, fix for each
    formula $\exists y\,\theta(y,\bar x)$ a Skolem map
    $F_\theta^{\mathcal M}$ on $M$: when a witness exists it chooses one,
    and otherwise it takes a fixed default value.  Starting with
    \[
       H_0=\widehat U\cup\{\varnothing^{\mathcal M}\},
       \qquad
       H_{m+1}=H_m\cup
       \{F_\theta^{\mathcal M}(\bar a):
          \bar a\in H_m^{<\omega},\ \theta\text{ a formula}\},       \tag{E2}
    \]
    put $K=\bigcup_{m<\omega}H_m$ and let $\mathcal K$ be the substructure
    of $\mathcal M$ induced on $K$.  There are only countably many formulas,
    and the set of finite tuples from a countable set is countable.  Hence,
    by induction on $m$, every $H_m$ is countable, and so is $K$.

    The set $K$ is closed under all the chosen Skolem maps.  If
    $\mathcal M\vDash\exists y\,\theta(y,\bar a)$ with $\bar a\in K$, all
    entries of $\bar a$ occur together in some $H_m$, and
    $F_\theta^{\mathcal M}(\bar a)\in H_{m+1}$ is a witness.  The
    Tarski--Vaught test therefore gives $\mathcal K\preccurlyeq\mathcal M$.
    This $K$ is the smallest set containing $\widehat U$ and the default
    value that is closed under the fixed Skolem maps.

    Moreover, the first downward-closure step in (E1) gives
    \[
     \operatorname{Ext}_{\mathcal K}(u)
    =\operatorname{Ext}_{\mathcal M}(u)\qquad(u\in U).       \tag{E3}
    \]
    Apply Theorem~\ref{thm:kaufmann-taller-extension} to $\mathcal K$.
    The resulting proper $\Sigma_{n+2}$-elementary taller* end extension
    $\mathcal N$ has the required agreement: for parameters from $U$, both
    $\mathcal M$ and $\mathcal N$ agree with $\mathcal K$ at level
    $\Sigma_{n+2}$.  Endness over $\mathcal K$, together with (E3), fixes
    every member of $U$, and the taller* witness for $K$ contains every
    member of $U$.  Since $K\subsetneq N$, also $U\subsetneq N$.  The
    full-Collection case uses all Skolem maps in exactly the same way and is
    otherwise identical.

Conversely, we again use induction on $n$ in order not to conceal the
bounded-quantifier calculation.  The standing assumption supplies
$\mathsf{Coll}(\Sigma_1)$ at the base, and the lower instances of (2)
follow from its $\Sigma_{n+2}$ instance.  By
Theorem~\ref{thm:collection-variants}, write the relevant
$\Pi_{n+1}$ formula as
\[
   \forall z\,\theta(x,y,z,\bar v),
   \qquad\theta\in\Sigma_n,
\]
and suppose
\[
 \mathcal M\vDash
 \forall x\in p\,\exists y\,\forall z\,
       \theta(x,y,z,\bar v).                               \tag{E4}
\]
Choose externally one witness $y_x$ for each
$x\in\operatorname{Ext}_{\mathcal M}(p)$.  Put $p,\bar v$, all these
$x,y_x$, and their downward closure into a countable set $U$.  In the
model supplied by (2), fixedness of $p$ prevents new members of $p$, while
partial elementary equivalence preserves each $\Pi_{n+1}$ witness
statement.  A set $c$ containing $U$ therefore gives
\[
 \mathcal N\vDash
 \exists q\,\forall x\in p\,\forall w\,
       \exists y\in q\,\forall z\in w\,
       \theta(x,y,z,\bar v).                               \tag{E5}
\]

For $n=0$, (E5) is visibly $\Sigma_2$.  If $n\geq1$, the induction
hypothesis gives $\mathcal M\vDash\mathsf{Coll}(\Sigma_n)$.  Its
$\Pi_{n+2}$ axiomatization transfers to $\mathcal N$ by the stated
$\Sigma_{n+2}$ equivalence.  In both models, the witness-collection
calculation (L3) converts the matrix after $q$ in (E5) to a
$\Pi_{n+1}$ formula.  Hence (E5) is $\Sigma_{n+2}$ and reflects to
$\mathcal M$.

Finally, if some $x\in p$ had no single $y\in q$ satisfying all
$z$-instances, then $\mathsf{Coll}(\Pi_n)$, equivalently
$\mathsf{Coll}(\Sigma_{n+1})$, would collect one counterexample for every
$y\in q$ into a set $w$.  This contradicts the corresponding instance of
(E5).  Thus $q$ collects genuine $\Pi_{n+1}$ witnesses in $\mathcal M$,
proving $\mathsf{Coll}(\Sigma_{n+2})$.  Running this induction through all
finite levels gives the full-Collection assertion.
\end{proof}

\subsection{Resolutions and stronger base theories}

The elementary role of a resolution is already visible before the
Keisler--Morley construction: in a sufficiently elementary proper end
extension, the first new stage covers the old universe.

\begin{theorem}\label{thm:resolution-makes-taller}
Let $n\geq0$.  Suppose that
$\mathcal M\preccurlyeq_{\Sigma_{n+2}}\mathcal N$ is a proper end
extension, $\mathcal M\vDash\mathsf{DB}_0$, and $\mathcal M$ has a
resolution of complexity at most $n+1$.  Assume that the ordinals of
$\mathcal N$ are linearly ordered by membership.  Then $\mathcal N$ is
taller than $\mathcal M$.
\end{theorem}
\begin{proof}
Let $r\leq n+1$ be the complexity of the resolution, and write
$R_\Sigma(x,\alpha)$ and $R_\Pi(x,\alpha)$ for its $\Sigma_r$ and
$\Pi_r$ definitions of stage membership.  We first isolate the formulas
that must transfer.  The covering assertion can be written
\[
  \forall x\,\exists\alpha\,\exists P\bigl(
     \Ord(\alpha)\wedge\rho_\Sigma(P,\alpha,v)\wedge x\in P
  \bigr).                                                   \tag{R1}
\]
For $r\geq1$ the two existential quantifiers merge into the leading
existential block of $\rho_\Sigma$, so (R1) is $\Pi_{r+1}$; for $r=0$ it
is $\Pi_2$.  In either case it is at most $\Pi_{n+2}$ and therefore holds
in $\mathcal N$.

For each old ordinal $\alpha$ and its unique old stage $P_\alpha$, the
following sentence holds in $\mathcal M$:
\[
 \forall x\bigl[
   (R_\Sigma(x,\alpha)\rightarrow x\in P_\alpha)
   \wedge
   (x\in P_\alpha\rightarrow R_\Pi(x,\alpha))
 \bigr].                                                    \tag{R2}
\]
For $r\geq1$, each implication in (R2) is $\Pi_r$ after prenexing and the
leading universal quantifier does not add an alternation; for $r=0$, (R2)
is $\Pi_1$.  Thus it is at most $\Pi_{n+1}$ and also transfers to
$\mathcal N$.  Similarly, monotonicity
has the transferred $\Pi_r$ form
\[
 \forall\beta\,\forall\gamma\,\forall x\bigl(
   \beta\in\gamma\wedge R_\Sigma(x,\beta)
      \longrightarrow R_\Pi(x,\gamma)
 \bigr).                                                    \tag{R3}
\]

Choose $c\in N\setminus M$.  By (R1), take an
$\mathcal N$-ordinal $\alpha_0$ and a stage $Q$ at $\alpha_0$ with
$c\in Q$.  If $\alpha_0\in M$, then (R2), with
$P_{\alpha_0}\in M$, gives $c\in P_{\alpha_0}$ in $\mathcal N$.
Endness would force $c\in M$, a contradiction.  Hence
$\alpha_0\notin M$.

Let $\beta$ be an old ordinal.  Linearity of the ordinal order in
$\mathcal N$ leaves, apart from equality, the alternatives
$\beta\in\alpha_0$ and $\alpha_0\in\beta$.  Equality is impossible, and
the latter alternative would make the new element $\alpha_0$ a member of
the old set $\beta$, contrary to endness.  Thus
$\beta\in\alpha_0$.

Finally fix $x\in M$.  The old covering assertion supplies an old
$\beta$ with $R_\Sigma(x,\beta)$.  Formula (R3) and
$\beta\in\alpha_0$ give $R_\Pi(x,\alpha_0)$ in $\mathcal N$, hence
$x\in Q$ by uniqueness of the stage at $\alpha_0$.  Therefore $Q$ contains
every element of $M$ and is the required taller witness.
\end{proof}

A resolution is also restrictive for large thin models.

\begin{theorem}\label{thm:aleph-one-resolution-collection}
Every $\aleph_1$-like model of $\mathsf{DB}_0$ which has a resolution
satisfies full Collection.
\end{theorem}
\begin{proof}
Let $\langle P_\alpha:\alpha\in\Ord^{\mathcal M}\rangle$ be the increasing
cover.  Its index order has no countable cofinal subset: otherwise
countably many sets, each with only countably many external members, would
cover $M$, contrary to $|M|=\aleph_1$.

Fix any formula $\varphi(x,y,\bar v)$ and suppose
\[
   \mathcal M\vDash
   \forall x\in p\,\exists y\,\varphi(x,y,\bar v).
\]
Externally choose a witness $y_x$ for each member of $p$.  This is a
countable family.  Choose $\alpha_x$ with $y_x\in P_{\alpha_x}$; the
countable set of indices has an upper bound $\alpha$.  Monotonicity puts
every $y_x$ into $P_\alpha$, which is a Collection bound.
\end{proof}

Thus the dead-end models from Theorem~\ref{thm:aleph-one-dead-end} may
behave quite differently depending on whether they carry a resolution.
We now prove the finite-strength resolved extension theorem.

For the rest of this subsection fix a standard $N\geq6$.  We first record
the set-existence and partition facts used in the construction.

\begin{lemma}\label{lem:B-axiomatization}
For $s\geq1$, the theory $\B_s$ has a $\Pi_{s+2}$ axiomatization.
Moreover, $\mathsf{Coll}_s(\Sigma_s)$ proves
$\mathsf{Sep}(\Sigma_s)$ and functional $\Sigma_s$-Replacement over
$\mathsf{DB}_0$.
\end{lemma}
\begin{proof}
By Theorem~\ref{thm:strong-collection-equivalence}, strong
$\Sigma_s$-Collection gives both ordinary $\Sigma_s$-Collection and
$\Sigma_s$-Separation.  Suppose
$\varphi(x,y,\bar p)\in\Sigma_s$ is functional on a set $a$.  Strong
Collection first gives $b$ such that every value of the relation on $a$
belongs to $b$.  The graph is then
\[
  G=\{\langle x,y\rangle\in a\times b:
          \varphi(x,y,\bar p)\}.                            \tag{B1}
\]
The defining formula in (B1) is $\Sigma_s$, so
$\Sigma_s$-Separation forms $G$.  Its range may either be obtained from
$G$ by the $\mathsf{DB}_0$ domain and coordinate operations or separated
from $b$ by
\[
   y\in b\ \wedge\ \exists x\in a\,
      \varphi(x,y,\bar p).
\]
The latter is a $\Delta_0(\Sigma_s)$ formula and is $\Sigma_s$ by
Lemma~\ref{lem:bounded-contraction}.  This proves functional
$\Sigma_s$-Replacement, including the existence of its graph.

For the axiomatization bound, the axioms of $\mathsf{DB}_0$ are
$\Pi_2$.  Infinity can be written as a $\Pi_2$ sentence by adding a dummy
leading universal quantifier to the usual existential assertion that an
inductive set exists.  Power Set has the form
\[
   \forall a\,\exists b\,\forall x
      \bigl(x\in b\leftrightarrow x\subseteq a\bigr),       \tag{B2}
\]
and is $\Pi_3$, since $x\subseteq a$ is bounded.

It remains to count strong Collection.  For $s=1$, formulas (A1)--(A3) in
the proof of Theorem~\ref{thm:collection-axiomatization} give a
$\Pi_3$ axiomatization.  For $s>1$, replace
$\Sigma_s$-Collection by $\Pi_{s-1}$-Collection.  If
$\psi(x,y)\in\Pi_{s-1}$, the nontrivial implication in a strong-Collection
bound has matrix
\[
   \forall y\,\neg\psi(x,y)
      \ \vee\ \exists y\in q\,\psi(x,y).                  \tag{B3}
\]
The first disjunct is $\Pi_s$ and the second contracts to
$\Pi_{s-1}$ over the lower Collection theory.  After the bounded
quantifier $x\in p$, the existential quantifier for $q$, and universal
closure, (B3) is $\Pi_{s+2}$.  The two inclusions defining a
$\Sigma_s$ separator are respectively $\Pi_s$ and $\Sigma_s$; the same
prenex calculation puts their universal closure in $\Pi_{s+2}$.  Hence
the strong-Collection part, and therefore all of $\B_s$, has a
$\Pi_{s+2}$ axiomatization.
\end{proof}

Replace any external sequence
$\langle\alpha_j:j<\omega\rangle$ by its finite running suprema.  For a
strict transitive resolution put
\[
   U_j=\operatorname{Ext}_{\mathcal M}(P_{\alpha_j}),
   \qquad U=\bigcup_{j<\omega}U_j.                         \tag{K1}
\]
Then $U_j\subseteq U_{j+1}$ and $U$ is externally transitive:
\[
   a\in U\ \wedge\ x\in^{\mathcal M}a
       \quad\Longrightarrow\quad x\in U.                 \tag{K2}
\]
Every finite tuple from $U$ lies in one $U_j$.  Also, whenever a set $b$
is formed in $\mathcal M$ by Separation, Power Set, strong Collection, or
finitely many functional Replacements, the resolution places $b$ in some
stage $P_\beta$; transitivity then gives $b\subseteq P_\beta$.  By taking a
successor of a finite supremum, $\beta$ can be placed above any prescribed
finite list of old bounds.  We call this observation \emph{stage capture}.

\begin{lemma}[Choice-free finite partition bounds]
\label{lem:choice-free-partition}
The theory $\B_6$ proves the following scheme, one instance for every
standard $m<\omega$.  Given a set $I$ of colors and a well-orderable
cardinal $\lambda$, there is an infinite well-orderable cardinal $\kappa$
such that every partition
\[
      [\kappa]^{m+1}=\bigcup_{i\in I}C_i
\]
has a homogeneous $Y\subseteq\kappa$ with $|Y|>\lambda$.
No form of Choice is required.
\end{lemma}
\begin{proof}
This is the finite-exponent theorem in the appendix to Section~4 of
\cite{Ke68}; we recall why its bounds exist at the stated strength.  For a
set $A$, let $\lVert A\rVert$ be the supremum of the initial ordinals
injecting into $A$.  For an infinite cardinal $\mu$, let $\mu^\circ$ be a
cardinal bounding the well-orderable cardinalities of unions of at most
$\mu$ sets, each of size at most $\mu$.  Power Set codes all relations and
functions on the relevant sets, Hartogs bounds their possible well-orders,
and functional Replacement takes the order types.  Thus $\mu^\circ$ exists
without well-ordering an arbitrary set.

More explicitly, first form the set of initial ordinals $\lambda$ for which
there are a sequence $\langle A_\xi:\xi<\mu\rangle$ and an injection
\[
  e:\lambda\longrightarrow\bigcup_{\xi<\mu}A_\xi,
  \qquad |A_\xi|\leq\mu\quad(\xi<\mu).                     \tag{P1}
\]
All candidate relations and functions in (P1) lie in Power Sets of fixed
products.  Hartogs' theorem bounds the possible domains $\lambda$, and
functional Replacement sends each coded well-order to its order type.
Choose $\mu^\circ$ strictly above the supremum of the resulting set of
initial ordinals.  It then has exactly the required union-bounding
property.

Starting with $\lVert I\rVert$, iterate finitely the operations
$\mu\mapsto\mu^\circ$, Power Set, and Hartogs, and finally choose $\kappa$
above the resulting bound and above $\lambda$.  For $m=0$, if every color
class were small, their union would be bounded by the defining property of
$\mu^\circ$, contradicting the size of $\kappa$.

For the induction step, recursively associate distinct coherent functions
$f_\alpha$ with $\alpha<\kappa$.  Initial segments of $f_\alpha$ are
represented by earlier $f_\beta$, and the next value records the original
color of a tuple obtained by adjoining $\alpha$.  If every domain were
below the previous bound, Power Set and $\mu^\circ$ would again bound the
number and union of all possible codes, contradicting $\kappa$.  Hence one
$f_\alpha$ has a sufficiently large domain.  Apply the induction
hypothesis to the induced coloring on the corresponding predecessor set;
coherence makes the resulting set homogeneous for the original coloring.

The actual nonempty color classes are well-orderable: map each to its least
member in the well-ordered domain $[\kappa]^{m+1}$.  This is the point at
which the argument avoids a choice function on $I$.

For completeness, the uniform syntactic bound is obtained as follows.
After the ambient Power Sets have been named, ``is a relation'', ``is an
injection'', and ``codes a well-order on a subset of $A$'' use only bounded
quantifiers over those Power Sets.  Separating the relevant codes and
forming their order types can conservatively be expressed by
$\Sigma_4$ formulas.  At the induction step, a candidate graph for the
sequence $\langle f_\alpha:\alpha<\gamma\rangle$ is again a member of a
fixed Power Set.  The assertion that every initial segment obeys the four
recursion clauses is $\Pi_5$ after the definitions of the code and the
least earlier representative are expanded.  Existence and uniqueness of
the next graph are therefore handled by functional $\Sigma_6$-Replacement.
These bounds also cover the finite iterations of $\mu\mapsto\mu^\circ$,
Power Set, and Hartogs.  Thus $\B_6$ proves every fixed finite-exponent
instance; no Set Foundation is used.
Adding $\lambda$ as a disjoint family of empty colors ensures that the
homogeneous output remains larger than any prescribed $\lambda$.
\end{proof}

The threshold $6$ is a deliberately conservative uniform bound for the
coding in this lemma.  It does not grow with term depth.  Improving this
fixed threshold is independent of the two-level loss in the output theory.

\subsubsection{Skolem envelopes and block homogenization}

For every formula
\[
   \exists y\,\theta(y,\bar x),\qquad \theta\in\Pi_{N-1}, \tag{K3}
\]
introduce a formal Skolem symbol $F_\theta$.  Lower-complexity formulas are
included by dummy quantifiers.  We do not choose single witnesses inside
$\mathcal M$.

\begin{lemma}[Bounded envelopes]\label{lem:bounded-envelopes}
Let $\mathcal M\vDash\B_N$ have a strict transitive resolution of
complexity at most $N$.  Given finitely many formulas (K3), finitely many
parameters, and an internal set $B$ of ordinal indices above their stages,
there is a higher stage $Q$ such that, for each relevant tuple
$\bar x\in\{P_\beta:\beta\in B\}^{<\omega}$,
\[
 E_\theta(\bar x)
   =\{y\in Q:\mathcal M\vDash\theta(y,\bar x)\}            \tag{K4}
\]
exists and is nonempty whenever the antecedent
$\exists y\,\theta(y,\bar x)$ holds.  The graphs of all these envelope
maps, and of their compositions through any fixed finite term depth, form
sets and remain $\Delta_N$-definable.
\end{lemma}
\begin{proof}
The stage relation has both a $\Sigma_N$ and a $\Pi_N$ definition and is
functional in the index.  Functional $\Sigma_N$-Replacement from
Lemma~\ref{lem:B-axiomatization} therefore forms
\[
    S_B=\{P_\beta:\beta\in B\}.                            \tag{K4a}
\]
Only finitely many arities occur in the data of the lemma.  Repeated
Cartesian Product forms the corresponding finite tuple powers of $S_B$;
there is no appeal here to a set of all finite sequences.

Let $D$ be the disjoint union of the finitely many sets of relevant pairs
$(\theta,\bar x)$.  The assertion
\[
   \exists y\,\theta(y,\bar x),
   \qquad\theta\in\Pi_{N-1},                               \tag{K4b}
\]
is $\Sigma_N$.  Strong $\Sigma_N$-Collection applied over $D$ supplies one
set $b$ containing a witness for every true instance (K4b).  Stage capture
then gives a transitive stage $Q$ containing $b$, all parameters, and all
sets used to code $D$.  Consequently every true instance has a witness in
$Q$.

For fixed $(\theta,\bar x)$, the envelope
\[
 E=\{y\in Q:\theta(y,\bar x)\}
\]
exists by $\Sigma_N$-Separation, since
$\Pi_{N-1}\subseteq\Sigma_N$ after adding a dummy existential quantifier.
Moreover, the assertion that $E$ is exactly this set can be written
\[
   E\subseteq Q\ \wedge\
   \forall y\in Q\bigl(y\in E\leftrightarrow
        \theta(y,\bar x)\bigr).                            \tag{K4c}
\]
After separating the two implications, bounded contraction puts (K4c) in
$\Delta_N$.  It defines $E$ uniquely, so functional
$\Sigma_N$-Replacement forms the indexed graph of all envelope values.

Now argue by induction on the fixed term depth.  Suppose graphs for the
immediate subterms have already been formed and are $\Delta_N$.  A possible
value of the composite term is specified by quantifying subterm values only
inside their envelopes and then applying one of the $\Delta_N$ graphs.
All newly introduced quantifiers are therefore bounded.  Applying
Lemma~\ref{lem:bounded-contraction} to both the $\Sigma_N$ and $\Pi_N$
definitions keeps the composite relation $\Delta_N$, and functional
Replacement forms its graph.  Thus no alternation is accumulated when the
term depth increases, proving all assertions of the lemma.
\end{proof}

Enumerate externally the $\Sigma_N$ Skolem schemes.  At block $s$ include
the first $s$ schemes, terms of depth and arity at most $s$, all parameter
and target substitutions from $U_s$, and all relevant subterms and
subtuples.  Although $U_s$ is an external collection of members, the
internal set $P_{\alpha_s}$ is a common domain for the entire block.

\begin{lemma}[Block homogenization]
\label{lem:bounded-block-homogenization}
Under the hypotheses of Lemma~\ref{lem:bounded-envelopes}, for each closed
parameter instance $\tau(\bar z,\bar u)$ and $a\in U$, choose externally
\[
 h(\tau,\bar u,a)\in
 \operatorname{Ext}_{\mathcal M}(a)\cup\{\star\}.         \tag{K5}
\]
These decisions are finitely realizable in the following strong sense.
Fix $q$.  Given finite sets of decisions, Skolem axioms, and old bounds,
there are
\[
   \beta_0\in\cdots\in\beta_{q-1}                         \tag{K6}
\]
above those bounds such that substitution of $P_{\beta_i}$ for the
distinct variables $z_i$ realizes
\[
\begin{cases}
 \tau(\bar P_\beta,\bar u)=h(\tau,\bar u,a),
       &h(\tau,\bar u,a)\ne\star,\\
 \tau(\bar P_\beta,\bar u)\notin a,
       &h(\tau,\bar u,a)=\star,
\end{cases}                                               \tag{K7}
\]
coherently for subterms and consistently with all the specified
$\Sigma_N$ Skolem witness axioms.
\end{lemma}
\begin{proof}
Proceed simultaneously through the blocks and through term depth.
Assume the earlier decisions are realizable on reservoirs above every
prescribed internal cardinal bound.  Start the next block with an ordinal
reservoir large enough for all the finitely many applications of
Lemma~\ref{lem:choice-free-partition} and above the old stage bounds.
Lemma~\ref{lem:bounded-envelopes} puts all possible values of the block's
terms in one stage $Q$.

It is important not to choose a value independently for each target.  Put
$A_s=P_{\alpha_s}$ and, for $v\in Q$, record the complete target profile
\[
  \pi_v:A_s\longrightarrow2,\qquad
  \pi_v(a)=1\ \longleftrightarrow\ v\in a.                \tag{K8}
\]
For an envelope $E\subseteq Q$ and a profile $\pi$, put
\[
  C_E(\pi)=\{v\in E:\pi_v=\pi\}.                          \tag{K9}
\]
Power Set forms all profiles and cells, and functional Replacement forms
their indexed family.  A \emph{coherent record} for an increasing tuple of
reservoir indices consists of the envelopes of all subterms, their
$\Pi_{N-1}$ witness relations, the composition relations, all earlier
restrictions, and all the cells (K9), simultaneously for every target in
$A_s$.  Possible records lie in the Power Set of one set-sized record
space.  The map taking a tuple to its nonempty set of coherent records
exists by Lemmas~\ref{lem:B-axiomatization} and
\ref{lem:bounded-envelopes}.  This set of records, together with the
records of all subtuples, is used as the color.

We record the complexity hidden in this coding.  A profile is a subset of
$A_s\times2$, and the assertion that it is the graph of the map in (K8) is
bounded once $v$ and $A_s$ are fixed.  Thus Power Set first supplies the
ambient profile space, and $\Delta_0$-Separation selects the profiles that
occur.  For fixed $E$ and $\pi$, membership in the cell (K9) is
\[
   v\in E\ \wedge\
   \forall a\in A_s\bigl(\langle a,1\rangle\in\pi
           \leftrightarrow v\in a\bigr),                   \tag{K9a}
\]
which is a bounded formula relative to the $\Delta_N$ envelope relation.
Hence the cell relation is $\Delta_N$, and functional
$\Sigma_N$-Replacement forms its graph.

The remaining fields of a record are tested as follows.  Membership in an
envelope is $\Delta_N$ by Lemma~\ref{lem:bounded-envelopes}; a Skolem
witness matrix is $\Pi_{N-1}$; composition quantifies possible subterm
values only inside their envelopes; and all earlier restrictions and
profile equalities are bounded.  Lemma~\ref{lem:bounded-contraction},
applied to the $\Sigma_N$ and $\Pi_N$ presentations, therefore makes
``$r$ is a coherent record for $\bar\beta$'' a $\Delta_N$ relation.  The
record space is a finite product of Power Sets of the sets just constructed.
Separation selects the nonempty coherent-record sets, and functional
$\Sigma_N$-Replacement maps each increasing tuple $\bar\beta$ to that set.
Thus the color map used below is a genuine set-valued $\Delta_N$ map; no
unbounded quantifier block is suppressed in passing from envelopes to
colors.

Apply Lemma~\ref{lem:choice-free-partition} finitely many times, from the
largest arity down.  Its size-preservation clause leaves a homogeneous
reservoir above any size required at the next block.  The common color is
a nonempty set of coherent records.  In the external metatheory, choose
one record and process it by term depth.  For a term, choose a nonempty
profile cell.  If its profile has value $1$ at some target $a$, every
member of the cell belongs to $a$, and homogeneity makes the same cell
available on every relevant tuple.  Choose one $h$ from it, use equality
for exactly the targets containing $h$, and use $\star$ for the others.
By (K2), this $h$ belongs to $U$ and therefore has a constant.  If the
profile is identically zero on the targets, values may vary with the
tuple, but every decision is the same negative $\star$ decision.  Since
subterms were processed first, composition and the Skolem witness
conditions remain coherent.

The induction invariant is not the existence of one fixed infinite
reservoir.  It is the assertion that every internal cardinal bound can be
exceeded by a reservoir realizing the decisions made so far.  The
size-preserving form of the partition lemma maintains this invariant.
After the external $\omega$-recursion, every finite collection of
requirements occurs in one block and is realized by a sufficiently long
finite sequence, proving (K6)--(K7).
\end{proof}

The use of Choice here is entirely external, at the same metatheoretic
level as Compactness.  Inside $\mathcal M$ only nonempty sets of coherent
records and set-valued envelopes are constructed.

\subsubsection{The finite-strength theorem}

\begin{theorem}[Finite-strength Keisler--Morley]
\label{thm:finite-keisler-morley}
Let $N\geq6$, let $\mathcal M\vDash\B_N$, and suppose that
$\mathcal M$ has a strict transitive resolution of complexity $r\leq N$.
For an external sequence
$\langle\alpha_j:j<\omega\rangle$ of ordinals of $\mathcal M$, put
\[
   U=\bigcup_{j<\omega}
     \operatorname{Ext}_{\mathcal M}(P_{\alpha_j}).
\]
For every external set-sized linear order $(X,<)$ there is a model
$\mathcal N$ such that:
\begin{enumerate}
    \item $(\mathcal N,u)_{u\in U}\equiv_{\Sigma_N}
          (\mathcal M,u)_{u\in U}$;
    \item every $u\in U$ is fixed;
    \item some $c\in N$ is transitive in $\mathcal N$ and contains every
          member of $U$;
    \item $(X,<)$ embeds isomorphically into
          $(N,\in^{\mathcal N})$; and
    \item $\mathcal N\vDash\B_{N-2}$.
\end{enumerate}
If the sequence of indices is cofinal in
$\Ord^{\mathcal M}$, then $U=M$ and, after identifying the named copy,
\[
  \mathcal M\preccurlyeq_{\Sigma_N}\mathcal N,
  \qquad
  \mathcal M\subsetneq_{\mathsf{end}}\mathcal N;
\]
this end extension is taller*.
\end{theorem}
\begin{proof}
Give every $u\in U$ a constant $c_u$.  The signed diagram
\[
\begin{split}
D_N(\mathcal M,U)={}&
\{\varphi(\bar c_{\bar u}):\varphi\in\Sigma_N,\
      \mathcal M\vDash\varphi(\bar u)\}\\
&{}\cup
\{\neg\varphi(\bar c_{\bar u}):\varphi\in\Sigma_N,\
      \mathcal M\vDash\neg\varphi(\bar u)\}
\end{split}\tag{K10}
\]
contains both truth values, the full atomic diagram, and all distinctness
literals.  Add the symbols $F_\theta$ from (K3) and their Skolem axioms
\[
 \forall\bar x\,
 \bigl(\exists y\,\theta(y,\bar x)
 \longrightarrow\theta(F_\theta(\bar x),\bar x)\bigr).    \tag{K11}
\]

Adjoin a new largest point $*$ to $X$, obtaining a nonempty order
$Y=X+1$, and add constants $d_y$ for $y\in Y$.  Let $T^*$ consist of
(K10)--(K11) and the following sentences:
\begin{enumerate}
    \item $d_y\in d_z$ exactly when $y<z$, together with
    $d_y\ne d_z$ for $y\ne z$ and $d_y\notin d_y$;
    \item $\Trans(d_y)$, $c_u\in d_y$, and $d_y\ne c_u$ for
    $y\in Y$ and $u\in U$;
    \item every instance of (K7), with an increasing tuple of $d_y$'s
    substituted for the variables and with the decision fixed by
    Lemma~\ref{lem:bounded-block-homogenization}.
\end{enumerate}
Repeated generators in the last clause are first identified and the
distinct ones put in increasing order.  Variable terms are included, so
the generators themselves cannot become new members of old targets.

We check finite satisfiability.  A finite subset mentions only finitely
many $d_y$'s, old constants, Skolem symbols, and decisions.  Use
Lemma~\ref{lem:bounded-block-homogenization} to choose corresponding
stages
\[
       P_{\beta_0},\ldots,P_{\beta_{q-1}}
\]
above the stages of all named old parameters.  Strictness gives
\[
 P_{\beta_i}\in P_{\beta_j}\quad\longleftrightarrow\quad i<j,
\]
and the stages are transitive, distinct, and not self-membered.  Their
height and transitivity put all finitely named members of $U$ into every
stage.  No selected stage equals a named $u\in U$: if
$P_{\beta_i}=u\in P_\gamma$ for an earlier $\gamma$, then (R) would give
$\beta_i\in\gamma$, contrary to $\gamma\in\beta_i$.

The coherent record supplies the finitely prescribed Skolem values and
decisions.  Outside the finitely specified tuples, external Choice extends
each partial assignment to a total function on $M$, choosing a witness
when the antecedent in (K11) holds and a fixed default otherwise.  This is
an expansion of the external structure $\mathcal M$, not an assertion that
$\mathcal M$ satisfies Choice.  Thus every finite subset of $T^*$ is
satisfiable.

By Compactness choose $\mathcal A\vDash T^*$.  Let $H$ be the closure of
\[
       U\cup\{d_y:y\in Y\}
\]
under all the $F_\theta$, and let $\mathcal N$ be the
$\mathcal L_\in$-reduct induced on $H$.  The bounded Tarski--Vaught
argument gives
\[
       \mathcal N\preccurlyeq_{\Sigma_N}\mathcal A.       \tag{K12}
\]
We give the bounded Tarski--Vaught induction explicitly.  Atomic formulas
have the same truth value in the induced substructure and in $\mathcal A$.
Suppose agreement has been proved for the proper subformulas of a bounded
formula and
\[
   \mathcal A\vDash\exists z\in b\,\delta(z,\bar a),
   \qquad b,\bar a\in H.                                  \tag{K12a}
\]
The matrix in (K12a) has a corresponding Skolem symbol.  Schemes of lower
complexity were included by dummy quantifiers.  Closure of $H$ under that
symbol therefore supplies a witness in $H$; the
induction hypothesis verifies its matrix in $\mathcal N$.  The reverse
direction is immediate from the induced relation and the induction
hypothesis.  Negation gives bounded universal quantifiers, so
$\mathcal N$ and $\mathcal A$ agree on all $\Delta_0$ formulas.

Now induct on $1\leq k\leq N$.  Let
$\exists y\,\theta(y,\bar a)$ be $\Sigma_k$, where
$\theta\in\Pi_{k-1}$ and $\bar a\in H$.  If it holds in $\mathcal A$, the
corresponding $F_\theta(\bar a)$ belongs to $H$ and satisfies $\theta$ by
(K11); agreement for the $\Pi_{k-1}$ matrix puts the witness in
$\mathcal N$.  Conversely, a witness in $H$ that works in $\mathcal N$
works in $\mathcal A$ by the same lower-level agreement.  Negating the
$\Sigma_k$ equivalence gives $\Pi_k$ agreement.  This simultaneous
induction proves (K12) without invoking Skolem functions above level $N$.

It remains to verify the five conclusions.  The signed diagram (K10) and
(K12) give
\[
  (\mathcal N,u)_{u\in U}\equiv_{\Sigma_N}
  (\mathcal M,u)_{u\in U}.                                \tag{K13}
\]
Since $\mathcal M\vDash\B_{N-2}$ and, by
Lemma~\ref{lem:B-axiomatization}, this theory has a
$\Pi_{(N-2)+2}=\Pi_N$ axiomatization, every one of its axioms is true in
$\mathcal M$.  Agreement for $\Sigma_N$ formulas in (K13), after negation,
also gives agreement for $\Pi_N$ formulas.  The parameter-free instances
therefore transfer one by one to $\mathcal N$, proving
$\mathcal N\vDash\B_{N-2}$.

Fix $a\in U$ and suppose $b\in^{\mathcal N}a$.  Every element of the hull
is a term value, so after ordering its distinct generators write
\[
 b=\tau(d_{y_0},\ldots,d_{y_{k-1}},\bar u),
 \qquad \bar u\in U^{<\omega}.
\]
The $\star$ alternative in (K7) contradicts $b\in a$.  Hence the equality
alternative applies and
\[
 b=h(\tau,\bar u,a)\in\operatorname{Ext}_{\mathcal M}(a).
\]
The reverse inclusion follows from the atomic diagram and (K2).  Thus
every $a\in U$ is fixed.

Take $c=d_*$.  The second group of axioms makes $c$ transitive and puts
every member of $U$ into it; transitivity is bounded and therefore remains
true in the induced hull.  The map $x\mapsto d_x$ embeds $X$, since the
positive and negative order literals and distinctness axioms give
\[
 x<x'\quad\longleftrightarrow\quad
 d_x\in^{\mathcal N}d_{x'}.
\]

Finally, if the index sequence is cofinal, the resolution's cover and
monotonicity give $U=M$.  Then (K13) is
$\mathcal M\preccurlyeq_{\Sigma_N}\mathcal N$, and fixedness of every old
set is precisely endness.  The extension is proper: if the taller*
witness $c$ were old, it would contain every old element.  The
$\Delta_0$ separator
\[
      r=\{x\in c:x\notin x\}
\]
would then be old and belong to $c$, giving
$r\in r\leftrightarrow r\notin r$.  This contradiction completes the
proof.
\end{proof}

The input loses no levels when term depth increases.  The output loses two
levels only because strong $\Sigma_s$-Collection is axiomatized at
$\Pi_{s+2}$ and the construction gives $\Sigma_N/\Pi_N$ agreement.  We do
not claim that this loss is optimal, nor that the output resolution is
preserved.

\subsubsection{The classical case}

The same construction, with no complexity truncation, gives the original
theorem over $\mathsf{ZF}$.

\begin{theorem}[Keisler--Morley]\label{thm:classical-keisler-morley}
Let $\mathcal M\vDash\mathsf{ZF}$ and let
$\langle\alpha_j:j<\omega\rangle$ be an external sequence of ordinals.
Put
\[
 U=\bigcup_{j<\omega}
   \operatorname{Ext}_{\mathcal M}
   \bigl((V_{\alpha_j})^{\mathcal M}\bigr).
\]
For every set-sized linear order $(X,<)$ there is
$\mathcal N\vDash\mathsf{ZF}$ such that:
\begin{enumerate}
    \item $(\mathcal N,u)_{u\in U}\equiv
          (\mathcal M,u)_{u\in U}$;
    \item every $u\in U$ is fixed;
    \item some transitive $c\in N$ contains every member of $U$; and
    \item $(X,<)$ embeds isomorphically into
          $(N,\in^{\mathcal N})$.
\end{enumerate}
\end{theorem}
\begin{proof}
Use the strict transitive resolution $P_\alpha=V_\alpha$ and run the proof
of Theorem~\ref{thm:finite-keisler-morley} with all formula schemes rather
than only those of complexity at most $N$.  Full Separation and Replacement
form the envelopes and record maps at every finite term depth, while Power
Set forms the color spaces.  The choice-free partition theorem from
\cite[Appendix to Section~4]{Ke68}, in the form used in
Lemma~\ref{lem:choice-free-partition}, supplies the homogeneous reservoirs.

Use the full signed elementary diagram and close the compactness model
under all Skolem functions.  The resulting hull is fully elementary in the
compactness model, so it satisfies $\mathsf{ZF}$ and has full elementary
agreement with the named copy of $\mathcal M$.  The term decisions prove
fixedness exactly as in the finite proof, the added top generator is the
transitive cover, and the remaining generators realize $X$.
\end{proof}

\begin{remark}
If the sequence is cofinal in $\Ord^{\mathcal M}$, then $U=M$, and
Theorem~\ref{thm:classical-keisler-morley} gives a proper elementary
taller* end extension of $\mathcal M$.  No Choice inside $\mathcal M$ is
used.  What the argument needs beyond bare $\mathsf{DB}_0$ is now explicit:
Power Set for color spaces, sufficient Collection/Replacement for envelope
and record maps, Infinity for the partition bounds, and a strict transitive
resolution for the homogeneous stage representatives.
\end{remark}

Consequently, an $\aleph_1$-like model with Replacement but without
Collection, of the kind studied in \cite{Gi16}, cannot carry a resolution
of the present kind: Theorem~\ref{thm:aleph-one-resolution-collection}
would turn the external countability of each old set into Collection bounds.

\section{Further Research}

The preceding results leave three groups of questions.

\begin{enumerate}
    \item Theorems~\ref{thm:db0-exact-cover} and
    \ref{cor:mutual-interpretability} show that $\mathsf{DB}_0$ and
    $\mathsf{PA}^-$ agree in two rather different respects.  It is natural
    to seek further shared extension, definability, or conservation
    properties.  At the same time, such results should locate the boundary
    between $\mathsf{DB}_0$ and Mathias's $\mathsf{GJ}_0$.  The latter is
    obtained by adding the rudimentary replacement operation $R_8$ to
    $\mathsf{DB}_0$ \cite[Remark~1.12 and
    Propositions~2.62--2.72]{Ma06}.  A useful model-theoretic dividing
    property should therefore detect formation of rudimentary images rather
    than merely the pairing, product, and separator operations already
    available in $\mathsf{DB}_0$.

    \item Corollary~\ref{cor:mutual-interpretability} does not identify the
    composites of the arithmetic and DAG interpretations.  Can mutual
    interpretability be strengthened to bi-interpretability, or to a
    natural intermediate notion in which one or both composites are
    definably isomorphic to canonical cuts or graph quotients?  Any such
    result must account for the noncanonical sharing and padding erased by
    bisimulation.

    \item Theorem~\ref{thm:finite-keisler-morley} leaves four quantitative
    and structural issues.  Can the output be improved from $\B_{N-2}$ to
    $\B_{N-1}$ or $\B_N$?  Can the conservative partition threshold
    $N\geq6$ be lowered by sharper Hartogs and recursion coding?  Can the
    exact membership chain (R) be replaced by a definable cofinal subsystem
    of stages while still realizing arbitrary linear orders?  Finally,
    under what additional elementarity or diagram conditions is the input
    resolution preserved in the output model?
\end{enumerate}

\bibliography{ref}

\end{document}